 \documentclass[letterpaper]{article}  

\usepackage{wrapfig}
\usepackage{graphics} 
\usepackage{epsfig} 
\usepackage{mathptmx} 
\usepackage{times} 
\usepackage{amsmath} 
\usepackage{amssymb}  
\usepackage{subfigure}
\usepackage{amsthm}
\usepackage{authblk}

\newtheorem{theorem}{Theorem}

\newtheorem{assumption}[theorem]{Assumption}
\usepackage{color}

\newtheorem{definition}[theorem]{Definition}

\newtheorem{proposition}[theorem]{Proposition}

{ 

\newtheorem{example}[theorem]{Example}
}

\begin{document}
\title{Normally hyperbolic surfaces based finite-time transient stability monitoring of power system dynamics}
\author[1]{Sambarta Dasgupta\thanks{dasgupta.sambarta@gmail.com}}
\author[2]{Umesh Vaidya\thanks{ugvaidya@iastate.edu}}
\affil[1]{Monsanto Company}
\affil[2]{ECPE Department, Iowa State University}
\maketitle
\begin{abstract}
In this paper, we develop a methodology for finite time rotor angle stability analysis using the theory of normal hyperbolic surfaces. The proposed method would bring new insights to the existing techniques, which are based on asymptotic analysis. For the finite time analysis we have adopted the Theory of normally hyperbolic surfaces. We have connected the repulsion rates of the normally hyperbolic surfaces, to the finite time stability. Also, we have characterized the region of stability over finite time window. The parallels have been drawn with the existing tools for asymptotic analysis. Also, we have proposed a model free method for online stability monitoring. 
\end{abstract}
\section{introduction} \label{sec_intro}
Transient stability analysis of rotor angle dynamics in power system is a classical problem with tremendous importance for safe and reliable operation of power system \cite{anderson2008power}.  Essentially the problem is to determine if the system, following a fault or disturbance, will reach a safe operating point.  Given the importance of this problem, there is extensive literature on the problem of transient stability analysis of power system. The classical results on this problem include Lyapunov function and energy function-based approaches. With the advancement in sensor technology in the form of Phasor Measurement Unit (PMU), the transient stability analysis of rotor angle dynamics has received renewed interest \cite{PMU_application1,PMU_application2,cutsem}.   The PMU provides time stamped real-time high resolution measurement data of power system state thereby providing an opportunity to develop real-time stability monitoring and control techniques \cite{PMU_phadake1,PMU_phadake2}. The goal of real-time transient stability monitoring is to determine if the system state following a fault or disturbance will reach the desired steady state based on measurement data over short time window of 2-6 sec \cite{jie_rotor}, \cite{dasgupta_volt}. The existing Lyapunov function and energy function-based methods are developed for asymptotic stability  analysis  and hence they are not suited for finite-time transient stability analysis. In this paper, we adopt techniques from geometric theory of dynamical systems for the development of short-term transient stability monitoring. In particular, we use the theory of normally hyperbolic invariant manifolds for the development of theoretical framework \cite{fenichel1971persistence,wiggins1994normally,haller2011variational}. We exploit the geometry of the phase space of the rotor angle dynamics in the development of the theoretical framework.  The theory of normally hyperbolic invariant manifolds allows us to characterize the rate of expansion and contraction for co-dimension one invariant manifolds. The co-dimension one manifold is said to be normally expanding (contracting) over the finite period of time if a normal vector to this manifold is expanding (contracting)  over the finite period of time.  Furthermore the extremum of expansion and contraction rate scalar field can be shown to be identified with the stability boundaries of the fixed point whose stability is under consideration. We show that the normal expansion and contraction rates can be used as the indicator of stability over the finite period of time. In particular, if the normal vector is expanding then the system behavior is deemed unstable and vice versa. The normal rates are also used for the computation of finite-time transient stability margin in real-time. The main contributions of this paper are as follows. We  provide mathematical rigorous foundation for finite-time transient stability monitoring of rotor angle dynamics. We propose computational framework for the computation of stability margin using real-time measurement data. \\
The main contribution of this paper is - 1. to adopt the theory of normally hyperbolic surfaces to address identification of finite time stability boundary, 2. to adopt normal expansion rate as a stability certificate for online stability monitoring, and 3. Adopt Lyapunov exponent for model-free stability monitoring for fast real time applications. \\
\section{Transient Stability Problem and phase space structure of power system}\label{section_problem_def}
Transient stability problem studies the stability of the rotor angle dynamics following a severe fault or disturbance. The transient stability time frame is typically  $3-5$ sec \cite{classification_stability}. For wide area swings the time interval of interest could be extended to $10-20$ sec. The discussion in this section follows closely from \cite{classification_stability}. The transient stability problem can be stated mathematically as follows \cite{chiang1987foundations}:
\begin{eqnarray}
\dot x&=&f_I(x),\;\;\;-\infty <t<t_F ,\label{eq1}\\
\dot x&=&f_F(x),\;\;\;t_F\leq t< t_P , \label{eq2}\\
\dot x &=&f(x),\;\;\;t_P\leq t<\infty . \label{eq3}
\end{eqnarray}
where $x\in \mathbb{R}^N$ is the state vector. Before the occurrence of fault, the system dynamics evolves with Eq. (\ref{eq1}). The fault is assumed to occur at time $t=t_F$ and the system undergoes structural change with the dynamics given by Eq. (\ref{eq2}). The duration of fault is assumed to be between time interval $[t_F,t_P]$. The fault is cleared at time $t=t_P$ and the state evolution is governed by the post fault dynamics i.e., Eq. (\ref{eq3}). Before the fault occurs the system is operating at some known stable equilibrium point $x=x_I$. At the end of the fault the state of the system is given by
\[x_P=\Phi_F(x_I,t_P-t_F) . \]
where $\Phi_F(x,t)$ is the solution of the system (\ref{eq2}) at time $t$ with initial condition $x$ at $t=0$. We assume that the post-fault system has a stable equilibrium point at $x=x_s$. The problem of transient stability is to determine if the post fault initial state $x_P$ with the post fault system dynamics, i.e., Eq. (\ref{eq3}), will converge to the equilibrium point $x_s$ i.e.,
\[\lim_{t\to \infty}\Phi(x_P,t)=x_s . \]
where $\Phi(x,t)$ is the solution of system (\ref{eq3}). The finite-time transient stability problem can then be defined as follows:
\begin{definition}[Finite-time transient stability] The finite-time transient stability problem is to determine if
\[\lim_{t\to \infty}\Phi(x_P,t)=x_s . \]
based on the state measurement data $x(t)$ over time interval $t\in[t_P,t_P+T]$  and  the model information about the post-fault dynamical system i.e., vector field $f$.
\end{definition}
The geometry of the state space will play an important role in the development of the framework for the finite-time transient stability monitoring. In the following section, we discuss the geometrical properties of the phase space dynamics of Eq. (\ref{eq3}).
\subsection{Phase space structure of swing dynamics}
Consider the system equation for the post-fault dynamics
\begin{eqnarray}
\dot x=f(x),\label{sys}
\end{eqnarray}
where, $f:\mathbb{R}^n \to \mathbb{R}^n$ is assumed to be $C^1$ vector field. We next introduce some definitions \cite{chiang1987foundations}.
\begin{definition}[Hyperbolic fixed point] A fixed point or equilibrium point of Eq. (\ref{sys}) is said to be hyperbolic if all the eigenvalues of the Jacobian obtained from linearizing the system Eq. (\ref{sys}) at the equilibrium point has non-zero real part.
\end{definition}
\begin{definition}[Type-$1$ saddle point] A fixed point is termed as stable, unstable, or saddle, if the real part of the eigenvalues of the Jacobian are respectively negative, positive, or both. The saddle point that have only one eigenvalue with positive real part is called as type-$1$ saddle point.
\end{definition}
\begin{definition}[Domain of attraction]
The domain of attraction for the post-fault equilibrium point $x_s$ is denoted by $A(x_s)$ and is defined as follows:
\[ A(x_s) := \{~ x ~ | ~ \underset{t \to \infty}{\lim} ~\Phi(x,t)  =  x_s ~ \}. \]
Also, we denote the boundary of the set $A(x_s)$ as $\delta A(x_s)$.
\end{definition}
\begin{definition}[Stable \& unstable manifold]
Let $x_i$ be the hyperbolic equilibrium point for system (\ref{sys}).
The stable manifold of the equilibrium point $x_i$ is denoted by $W_s(x_i)$ and is defined as,
\[ W_s(x_i) := \{ ~ x ~ | ~ \underset{t \to \infty}{\lim} ~\Phi(x,t)  =  x_i ~ \}. \]
Similarly the unstable manifold of the equilibrium point $x_i$ is denoted by $W_u(x_i)$ and defined as,
\[ W_u(x_i) := \{ ~ x ~ | ~ \underset{t \to -  \infty}{\lim} ~\Phi(x,t)  =  x_i ~ \}. \]
\end{definition}
\begin{definition}[Transversal intersection] Consider two manifolds $A$ and $B$ in $\mathbb{R}^n$. We say that the intersection of $A$ and $B$ satisfies the transversality condition if at the point of intersection $x$, the tangent spaces of $A$ and $B$ at $x$ span the $\mathbb{R}^n$ i.e.,
\[T_x(A)+T_x(B)=\mathbb{R}^n . \]
or the manifold do not intersect at all. We use $T_x(A)$ to denote the tangent space of $A$ at point $x$.
\end{definition}
We now make following assumptions on system equation (\ref{sys}) as could be found in \cite{chiang1987foundations}.
\begin{assumption} \label{ass_sys}\[\]
\begin{itemize}
\vspace{-0.25in}
\item A1. All equilibrium points on the stability boundary of (\ref{sys}) are hyperbolic.
\item A2. The intersection of $W^s(x_i)$ and $W^u(x_j)$ satisfies the tranversality condition, for all the equilibrium points $x_i$, $x_j$ on the stability boundary.
\item A3. There exists a $C^1$ function $V: \mathbb{R}^n\to \mathbb{R}$ for (\ref{sys}) such that
\begin{itemize}
\item $\dot V(\Phi(x,t))\leq 0$.
\item If $x$ is not an equilibrium point then the set $\{t\in \mathbb{R}: \dot V(\Phi(x,t))=0\}$ has measure zero in $\mathbb{R}$.
\item $V(\Phi(x,t))$ is bounded implies that $\Phi(x,t)$ is bounded.
\end{itemize}
\end{itemize}
\end{assumption}

Assumptions A1 and A2 on the hyperbolicity and transversitilty are generic property of dynamical systems and hence true for almost all dynamical system. If any dynamical system does not satisfy assumptions A1 and A2 then any small perturbations of such dynamical system will satisfy assumptions A1 and A2. Assumption A3 guarantee that all the trajectories of the system will go to infinity or converge to one of the equilibrium point and hence the possibilities of limit cycle oscillations and chaotic motion is ruled out.

It is important to emphasize that both the classical power swing model and the structure preserving model of power system satisfy Assumption A1-A3. These assumptions have important consequences on the topological property of the boundary of domain of attraction of fixed point $x_s$. These consequence are presented in the form of following two fundamental theorems \cite{varaiya1985direct}.
\begin{theorem}
For a dynamical system (\ref{sys}) satisfying assumptions A1-A3, $x_i$ is an unstable equilibrium point on the boundary of domain of attraction i.e., $\delta A(x_s)$ of the stable equilibrium point $x_s$ if and only if $W^u(x_i)\cap A(x_s)\neq \emptyset$.
\end{theorem}
\begin{theorem} For a dynamical system (\ref{sys}) satisfying assumptions A1-A3, let $x_i$ for $i=1,2,\ldots$ be the unstable equilibrium points on the stability boundary $\delta A(x_s)$ of the stable equilibrium point $x_s$, then
\[\delta A(x_s)=\cup_{x_i\in E\cap \delta A} W^s(x_i),\]
where, $E$ is the set of all equilibrium points for (\ref{sys}).
\end{theorem}
The geometry of the boundary of domain of attraction for the stable fixed point $x_s$ for the post-fault dynamics play an important role in the development of theoretical foundation for finite-time transient stability monitoring. We discuss the theory necessary for the development of this framework in the following section.
\begin{figure}
    \centering
    \includegraphics[width=0.65\textwidth]{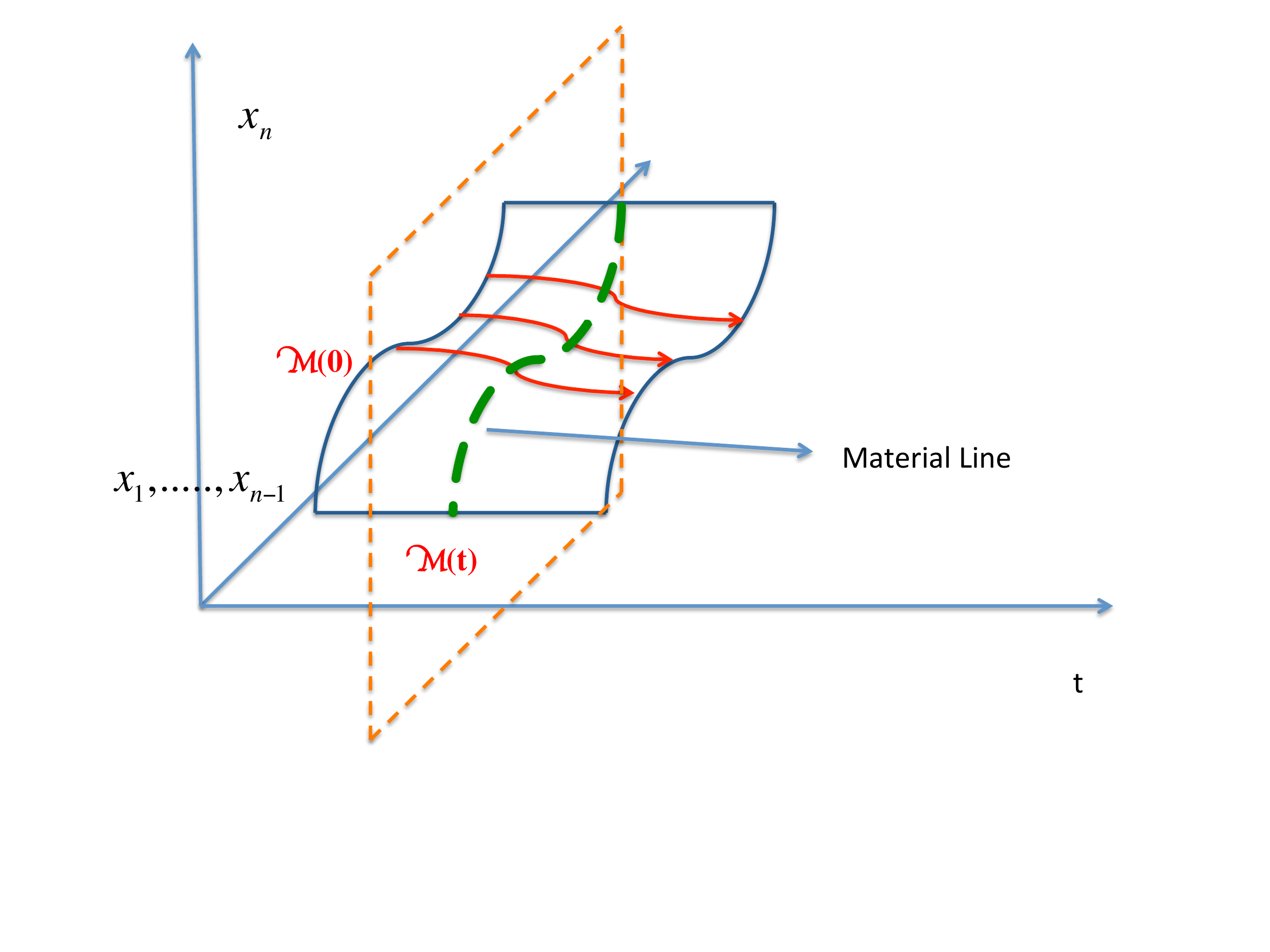}
  \caption{ Schematic of the material surface. }
 \label{fig_mat_surf_sche}
\end{figure}
\section{Normally Hyperbolic Invariant Manifold Theory}\label{section_normal_hyperbole}
In this section, we introduce some preliminaries for the theory of the normally hyperbolic invariant manifolds. For more details, we refer the readers to \cite{haller2011variational}. The theory of normal hyperbolicity over finite-time is developed for more general time varying system. However given our interest in system Eq. (\ref{sys}) for the post-fault dynamics which is time invariant, we restrict our discussion to time invariant vector field (\ref{sys}). We start with the following definition of material surface \cite{haller2011variational}.
\begin{definition} [Material surface \cite{haller2011variational} ]
 A material surface $ \mathcal M( t) $ is the $ t = $ const. slice of an invariant manifold  in the extended phase space $ X \times [\alpha, \beta] $, generated by the advection of an $n - 1$  dimensional surface of initial conditions $\mathcal{M} (0)$ by the flow map $\Phi ( x , t )$,
\begin{align*}
\mathcal{M} ( t ):= \Phi (\mathcal{M} (0), t),\;\;\dim \left ( \mathcal{M} (0) \right ) =  n-1, \;\;t\in[0,T].
\end{align*}
The schematic of the material surface is shown in Fig. \ref{fig_mat_surf_sche}.
\end{definition}
We want to express the attraction and repulsion property of this material surface over the time interval $[0,t]$. To this end we consider an arbitrary point $x_0\in {\cal M}(x_0)$ and $(n-1)$ dimensional tangent space $T_{x_0}{\cal M}(0)$ of ${\cal M}(0)$ and one dimensional normal space $N_{x_0} {\cal M}(0)$. The tangent vector is carried forward by along the trajectory $\Phi(x_0,t)$ by the linearized flow map given by $\nabla \Phi(x_0,t)$  into the tangent space
\[T_{\Phi(x_0,t)}=\nabla \Phi(x_0,t)T_{x_0}{\cal M}(0) . \]
By the invariance property of the manifold ${\cal M}(0)$, the tangent vector at point $x_0$ of ${\cal M}(0)$ are propagated to the tangent vector at point $\Phi(x_0,t)$ under the linearized flow  $\nabla \Phi(x_0,t)$. However this is not the case with the normal vector at point $x_0$ to manifold ${\cal M}(0)$. A unit normal vector $n_{x_0}\in N_{x_0}{\cal M}(0)$ will not necessarily be mapped into the normal space $N_{\Phi(x_0,t)}{\cal M}(t)$.
 \begin{definition}[Normal repulsion rate \cite{haller2011variational}]\label{rho_def}
 Let $n_{\Phi(x_0,t)}\in N_{\Phi(x_0,t)}{\cal M}(t)$ denotes the family of smoothly varying family of unit normal vector along the flow $\Phi(x_0,t)$. The growth of perturbation in the direction normal to ${\cal M}(t)$ is given by the repulsion rate and is given by
 \begin{eqnarray}
 \rho(x_0,n_0,t)=\left<n_{\Phi(x_0,t)},\nabla \Phi(x_0,t)n_0\right> .
 \end{eqnarray}
 \end{definition}
 If $\rho(x_0,n_0,t)>1$, then the normal perturbation to ${\cal M}(0)$ grows over the time interval $[0,t]$. Similarly, if $\rho(x_0,n_0,t)<1$, then the normal perturbations to ${\cal M}(0)$ decreases.
\begin{definition}[Normal repulsion ratio \cite{haller2011variational}]
The normal repulsion ratio is a measure of the ratio between the normal and tangential growth rate along ${\cal M}(t)$ and is defined as follows:
\begin{eqnarray}
v(x_0,n_0,t)=\min_{|e_0|=1, e_0\in T_{x_0}{\cal M}(0)}\frac{\left<n_{\Phi(x_0,t)}, \nabla \Phi(x_0,t)n_0\right>}{|\nabla \Phi(x_0,t)e_0|} .
\end{eqnarray}
\end{definition}
If $v(x_0,n_0,t)>1$, then the normal growth along ${\cal M}(t)$ dominates the largest tangential growth rate along ${\cal M}(t)$. From the point of view of computation, both the normal rate and ratio can be computed using the following formulas \cite{haller2011variational}.
\begin{eqnarray} \label{rho_formula}
 \rho(x_0,n_0,t)=\frac{1}{\sqrt{\left<n_0, [C_t(x_0)]^{-1}n_0\right>}} ,
\end{eqnarray}
where, $C_t(x_0):=\nabla \Phi(x_0,t)^{*}\nabla \Phi(x_0,t)$.
\begin{eqnarray}
v(x_0,n_0,t)=\min_{|e_0|=1,e_0\in T_{x_0}{\cal M}(0)} \frac{\rho(x_0,n_0,t)}{\sqrt{\left<n_0, [C_t(x_0)]^{-1}n_0\right>}} .
\end{eqnarray}
With the aid of normal repulsion rate and ratio, we define finite-time repelling (attracting) material surfaces.
\begin{definition} [Normally repelling material surface\cite{haller2011variational}]
A material surface $\mathcal{M} (t) $ is normally repelling over time interval $[0,T]$ if there  exist numbers $\alpha, \beta > 0$, such that, for all points $x_0 \in \mathcal{M}(0)$, and unit normal vector $n_0 \in N_{x_0}  \mathcal{M} (0) $, we have \cite{haller2011variational},
\begin{align*}
\rho ( x_0, n_0, t) & \ge e^{\alpha t} , \\
v ( x_0, n_0, t)& \ge e^{\beta t} .
\end{align*}
A material surface $\mathcal{M}(t)$ is normally attracting over time interval $[0,T]$ if it is normally repelling over $[0,T]$ in backward time. We call ${\cal M}(t)$ hyperbolic over $[0,T]$ if it is normally attracting or normally repelling over $[0,T]$.
\end{definition}
All repelling or attracting material surfaces are not equally repelling or attracting. Some material surfaces are more repelling or attracting than others. Material surfaces that are maximally attracting or repelling occupy special place and name in the theory of normally hyperbolic invariant manifolds. They are called as Lagrangian coherent structures. We have following definitions in this direction.
\begin{definition} [Lagrangian Coherent Surface \cite{haller2011variational}] \label{lcs_def}
Assume that ${\cal M}(t)$ is a normally repelling
material surface over $t\in [0,T]$. We call ${\cal M}(t)$ a repelling
Lagrangian Coherent Surface (LCS) over $[0, T]$ if its normal repulsion rate admits
a point-wise non-degenerate maximum along ${\cal M}(0)$ among
all locally $C^1$-close material surfaces. We call ${\cal M}(t)$ an
attracting LCS over $[0,T]$ if it is a repelling LCS over $[0, T]$
in backward time. We call ${\cal M}(t)$ a hyperbolic LCS over $[0,T]$ if it is repelling or attracting LCS over $[0,T]$.
\end{definition}
\begin{figure}[h]
\begin{center}
\mbox{
\hspace{-0.2in}
\subfigure[]{\scalebox{.8}{\includegraphics[width=3.0 in, height=2.0 in]{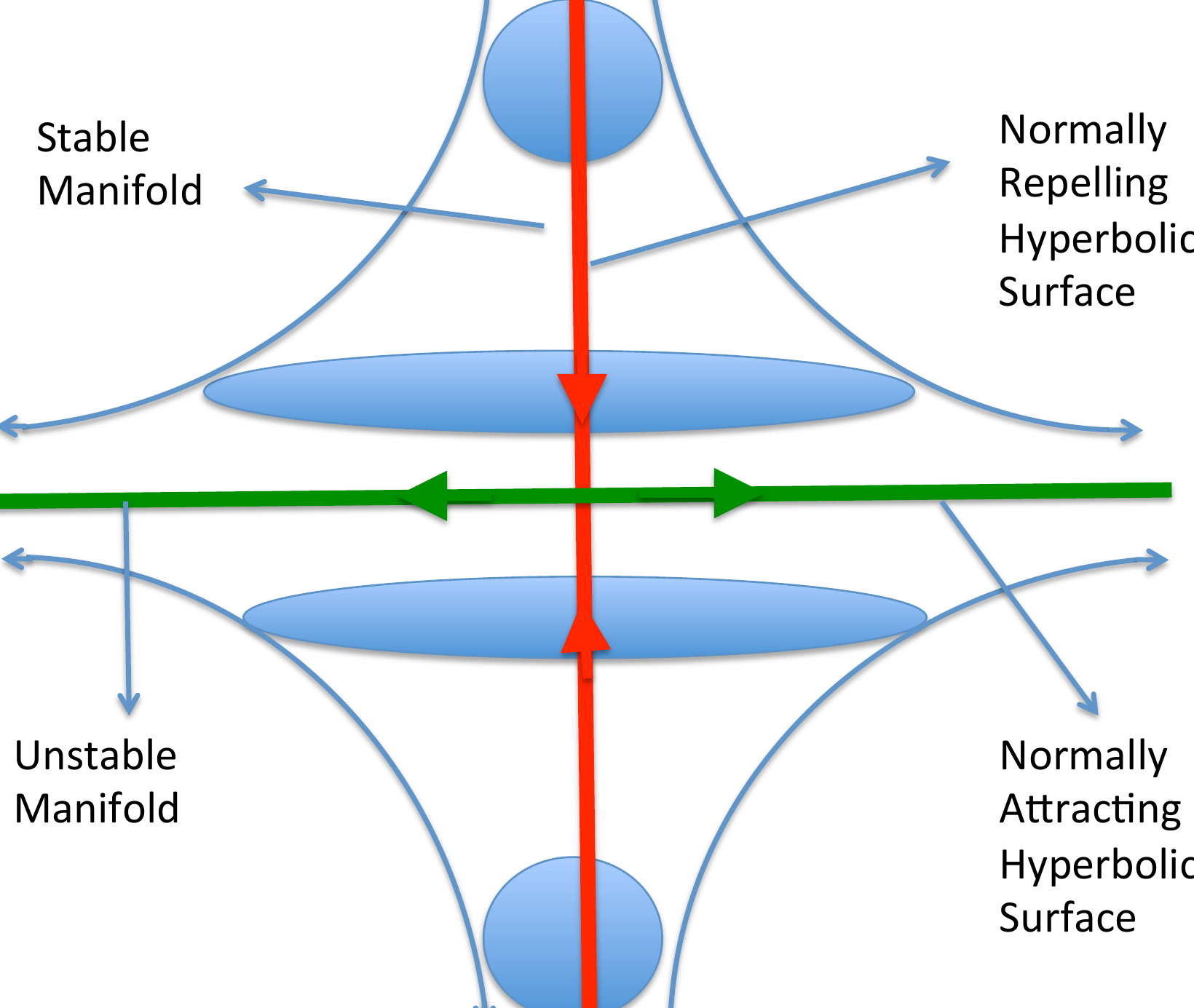}}}
\hspace{-0.1in}
\subfigure[]{\scalebox{.26}{\includegraphics[width=8.0 in, height=6.00 in]{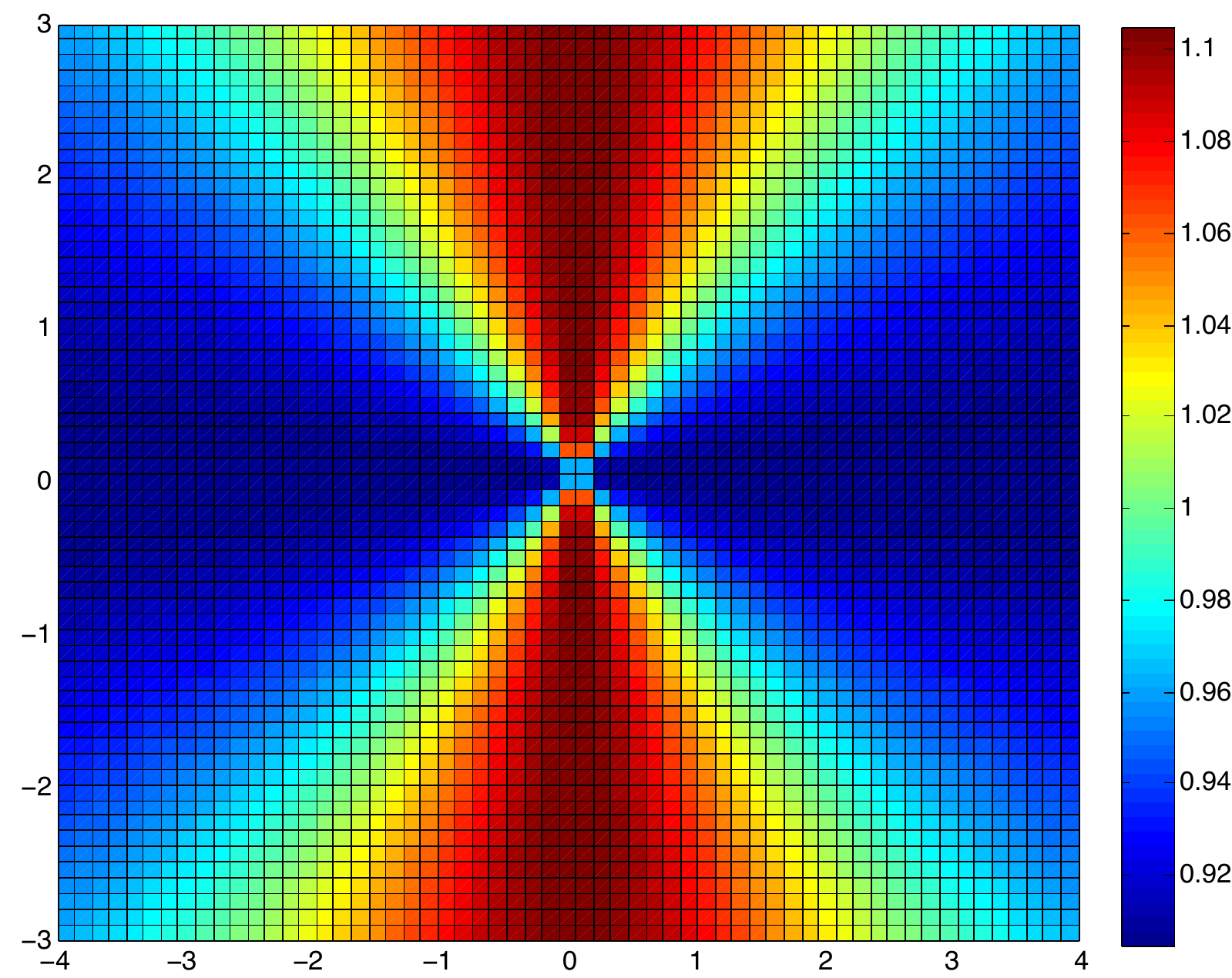}}} }
\caption{ (a) Normally repelling and attracting surfaces for saddle point, (b) Normal repulsion rate computed after $1$ sec. }
\label{fig_saddle_sche}
\end{center}
\end{figure}
With the aid of the following example we demonstrate the concept of the attracting and repelling surfaces.
\begin{example}
We consider the system,
\begin{align*}
\dot x  = x,\;\;\;
\dot y  = -y.
\end{align*}
It can be observed that the $y$ and $x$  axes form the stable and unstable manifolds of the saddle point.  As we can observe from Fig. \ref{fig_saddle_sche} (a) , the stable manifold $y$ axis repels the nearby trajectories and forms a repelling material surface. On the other hand the unstable manifold $x$-axis forms the attracting material surface.  This we verify by computing the normal expansion rates. The normal expansion rate is more than $1$ over a surrounding region of $y$ axis. The  $y$  axis form the Lagragian Coherent Structure for the system, as the normal expansion rate achieves a local maxima along the $y$ axis.
The normal expansion rate turns out to be.
\begin{equation*}
 \rho ( x_0, n_0, t ) = \frac{1}{ \sqrt {  n_0^T  \left ( \left( \nabla  \Phi(x_0, t)) \right)^T \nabla  \Phi(x_0, t)  \right ) ^{-1} n_0 }}  ~,
 \end{equation*}
 Where,
 \begin{align*}
 \left ( \left( \nabla  \Phi(x_0, t)) \right)^T \nabla  \Phi(x_0, t)  \right ) =  \begin{bmatrix}
 e^{2t} & 0 \\
 0 &  e^{ - 2t}
 \end{bmatrix}.
 \end{align*}
 For this system the surfaces $(0,y)^T$, and $(x,0)^T$ forms the repelling and attracting material surfaces, as we can observe in Fig. \ref {fig_saddle_sche} (b). Also, it can be observed there is a small patch of region around the $y$ - axis, which has normal expansion rate more than $1$. This means those material lines are repelling for the finite time of our interest. 
\end{example}
\subsection{Main Results} \label{section_mainresults}
We would relate the stability boundary to the normally hyperbolic surfaces with the aid of the following two Theorems.  These two theorems show that in finite time the normal expansion rate can prescribe a region around the stable manifold of a type -$1$ saddle point as repelling. Theorem \ref{norm_hyp_1} shows that the normal expansion rate can be used to characterize the finite time stability boundary, and Theorem \ref{norm_hyp_2} implicates that the normal expansion rate can be used as a stability certificate for online monitoring.
\begin{theorem} \label{norm_hyp_1}
For any finite time interval $T > 0$ and $x_0 \in W_s (x^e)$, there exists an $\epsilon > 0$ such that, following condition is satisfied,   \[ \rho ( x , \eta_0, T) > 1 , ~~ \forall x \in \mathcal{B}_{\epsilon} (x_0), \]
where,  $ \mathcal{B}_{\epsilon} (x) $ represents an open ball of radius $\epsilon$ around $x_0$, and $W_s ( x^e )$ denotes the stable manifold of  a type-$1$ saddle point $x^e$, which is located at the boundary of the domain of attraction $A(x^s)$.
\end{theorem}
$~~~$Next, we outlline Theorem \ref{norm_hyp_2}, which can specify the finite time stability region contained in the domain of attraction in terms of the normal expansion rates.
\begin{theorem} \label{norm_hyp_2}
For any finite time interval  $T > 0$ there exists an $\epsilon > 0$ such that,  following condition is satisfied,  \[ \rho ( x , \eta_0, - T ) > 1, ~ ~  \forall x \in  \mathcal{B}_{\epsilon} (x^s), \]
where,  $ \mathcal{B}_{\epsilon} (x^s) $  represents an open ball of radius $\epsilon$ around the stable fixed point $x^s$.
\end{theorem}

Theorem \ref{norm_hyp_1} provides a region for finite time stability based on the repulsion rate. Theorem \ref{norm_hyp_1}  also can be rephrased as following - for an arbitrary finite time interval $T>0$, there exists a $\epsilon > 0$ such that, a material surface at a  maximum distance of $\epsilon$ from the  the stable manifold of any type - $1$ saddle point , forms a finite time repelling hyperbolic surface.  Hence, the normal repulsion rate provides a rigorous way to demarcate the finite time stability boundary.  The classical methods were capable of demarcating the stability boundaries in terms of the asymptotic dynamics. But this may not be adequate for finite time stability analysis.  For finite time stability monitoring, we compute normal repelling rate of the trajectory and ascertain the stability. \\
$~~~$Theorem \ref{norm_hyp_2} describes the material surfaces that form the finite time repelling material surfaces.  for an arbitrary finite time interval $T>0$, there exists a $\epsilon > 0$ such that, a material surface inside the domain of attraction $ A(x^s)$ of the stable fixed point and at minimum distance of $\epsilon$ away from the all the stable manifolds of the type - $1$ saddle point  forms a finite time attracting hyperbolic surface. Also, the implication of Theorem \ref{norm_hyp_2} is that the normal repulsion rate can be used for a finite time stability certificate for online stability monitoring. \\
$~~~$  For both Theorem \ref{norm_hyp_1}, and  \ref{norm_hyp_2}, the parameter $\epsilon$ is a function of the time interval of interest, i.e. $T$.  The normal expansion rate can be used to obtain the finite time stability boundaries, and also can be used for online stability monitoring. It is to be noted that, the $\epsilon$ in Theorem \ref{norm_hyp_1}, decreases with increase in $T$, and $\epsilon$ in Theorem \ref{norm_hyp_2} increases with increase in $T$. The normal repulsion rate can be used as an estimate of the stability margin. Theorem \ref{norm_hyp_2} describes the material surfaces that are finite time repelling. Next, we would discuss the new insights that this approach brings to the existing method.

\subsection{ Margin of Stability from Normal Expansion Rate}
The existing techniques of stability analysis are based on energy function based methods. Energy function is a non-negative function defined over the state space, which guarantees asymptotic boundedness of the system trajectories. The energy function is bounded and decays along a trajectory, as a result of that the trajectories of the post fault system trajectories are bounded \cite{varaiya1985direct}. The Assumptions \ref{ass_sys} would imply that the trajectories inside domain of attraction would converge to the stable fixed point. This rules out possibility of limit cycle or other type of behavior inside the domain of attraction. A critical value of the energy function is given by the energy function value, evaluated at the type -$1$ saddle point. Energy values less than the critical level results in asymptotic stability. However, the energy functions are inadequate to specify the region of stability over a finite time. Also, energy function based techniques have extensively used for real time stability monitoring, by augmenting it with power flow type approaches \cite{5356818,726889}. We demonstrate that our approach is also suitable for asymptotic stability monitoring. In order to achieve this goal, we define the following function over the trajectory. \\ We first define a finite time stability margin inverse of  accumulation of repulsion rates over the trajectory. Stability margin corresponding to the point  $x_0 \in \mathcal{M}(0)$ over a time window $T$ is denoted as $\frac{1}{\gamma (x_0, \eta_0, T)}$ where,
\[ \gamma ( x_0, \eta_0, T ) := \int_{0}^{T} \rho (x_0, \eta_0, \tau) d \tau.\]
It is to be noted that $\eta_0 \in \mathcal{N}_{x_0} \left ( \mathcal{M} (0) \right )$ is unit normal vector. By taking $T \to \infty$, ~ $\gamma ( x_0, \eta_0, T ) $ reaches a constant value if the system is stable. This will result in a positive value of the margin. Larger the margin, more the system is stable. On the other hand if the system is unstable, the margin will go to $0$. \\ A non-negative function corresponding to the point  $x_0 \in \mathcal{M}(0)$ is denoted as $V_{\rho}(x_0)$ and is defined as,
\[ V_{\rho}(x_0) := \underset{T \to \infty}{\text{Lim}} ~ \gamma (x_0, \eta_0, T) = \int_{0}^{\infty} \rho ( x_0, \eta_0, \tau) d \tau. \]
\begin{proposition}
$V_{\rho}(x(t))$ is a non-negative bounded function, which is decreasing over the trajectories with initial conditions inside the domain of attraction \cite{varaiya1985direct}. Our approach can compute the margin also in finite time.
\end{proposition}
This guarantees the boundedness of the trajectories asymptotically inside the domain of attraction, similar to that of \cite{varaiya1985direct}. 
\section{Finite Time Lyapunov Exponent}
\begin{figure} [h]
    \centering
    \includegraphics[width=0.55\textwidth]{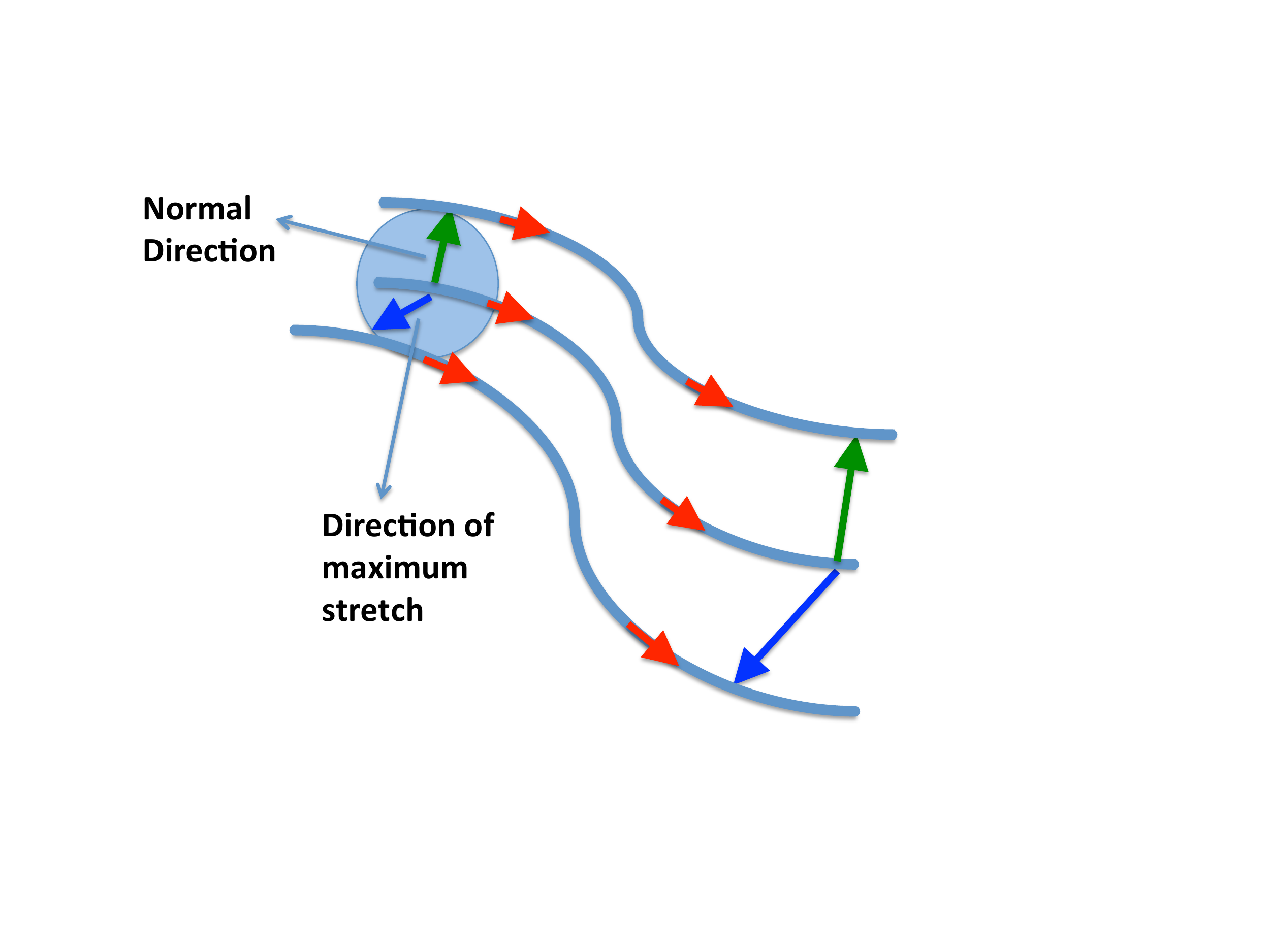}
  \caption{ The computational Scheme for normal expansion rate and Finite Time Lyapunov Exponent. }
 \label{fig_phase_portrait}
\end{figure}
Previously, we have shown that the normal expansion rates can detect the stability boundaries accurately, as stability boundaries form normally repelling hyperbolic structure. Also, we have demonstrated the normal expansion rate can be used as a stability certificate. Normal expansion rate computes the expansion in the normal vector along a trajectory. In this section, we would introduce Lyapunov exponent (LE), which is also a similar estimate of rate of expansion like normal repulsion rate. Finite Time Lyapunov Exponent (FTLE) computes the maximum local stretching for an initial condition \cite{hadjighasem_journal}, and it is also capable of identifying the normally repelling hyperbolic surfaces (i.e. stability boundaries) under some technical condition given in \cite{haller2011variational}. 
First we introduce the definition of Lyapunov exponent, and successively describe the finite time Lyapunov exponent (FTLE). Mathematical definition of maximum or principal Lyapunov exponent \cite{Katok} in asymptotic sense is is as follows,
\begin{definition}
Consider a continuous time dynamical system $\dot x=f(x)$, with $x\in X\subset \mathbb{R}^N$. Let $\Phi( x , t )$ be the solution of the differential equation. Define the following limiting matrix,
\begin{eqnarray}
\bar \Lambda(x) : =\lim_{t\to \infty}\left[\frac{\partial \Phi (x , t ) }{\partial x}^{T}\frac{\partial \Phi ( x, t )}{\partial x}\right]^{\frac{1}{2t}}. \label{limiting_matrix}
\end{eqnarray}
Let $\bar \Lambda_i(x)$ be the eigenvalues of the limiting matrix $\bar \Lambda(x)$. The Lyapunov exponents $\bar \lambda_i(x)$  are defined as
\begin{eqnarray}
\bar \lambda_i(x) : =\log \bar \Lambda_i(x).
\end{eqnarray}
Let $\bar \lambda_1(x)\geq \bar \lambda_2(x)\cdots, \geq \bar \lambda_N(x) $, then $\bar \lambda_1(x)$ is called the maximum Lyapunov exponent.
\end{definition}
Using results from Multiplicative Ergodic Theorem, it is known that the limit in Eq. (\ref{limiting_matrix})  is well defined  \cite{RevModPhys.57.617}. Furthermore, the limit in Eq. (\ref{limiting_matrix}) is independent of the initial condition, $x$, under the assumption of unique ergodicity of the system.  Lyapunov exponents can be thought of as the generalization of eigenvalues from linear systems to nonlinear systems in asymptotic sense. The Finite Time Lyapunov Exponent (FTLE) is defined as  \cite {Shadden_ftle}, 
\begin{definition}
The Finite-Time Lyapunov Exponent (FTLE) $\sigma^{\tau}(x_0)$ at an initial point $x_0$ for  time interval $\tau$ is defined as,
\begin{equation*} 
\sigma^{\tau}(x_0) := \frac{1}{\tau} \ln \parallel \frac{\partial \Phi}{\partial x} (x_0, \tau) \parallel.
\end{equation*} 
\end{definition}
 Figure \ref{fig_phase_portrait} shows conceptually the computation of the normal expansion rate and FTLE. An small ball around an initial condition is taken, and the initial conditions inside the ball are considered. After finite time the evolving distance between any pair of initial condition inside the ball would be studied. The green vector shows the evolving normal direction, and will thus be used to compute the normal expansion rate. On the other hand, the green vector gives the direction, which gets maximum stretched in the in the finite interval, and thus be used to compute the FTLE. Using the multiplicative ergodic theorem, it can be shown asymptotically the Lyapunov exponent becomes independent of the initial conditions \cite{RevModPhys.57.617}.    
 \section{Relation Between FTLE and Hyperbolic Surfaces}
\begin{figure} [h]
    \centering
    \includegraphics[width=0.55\textwidth]{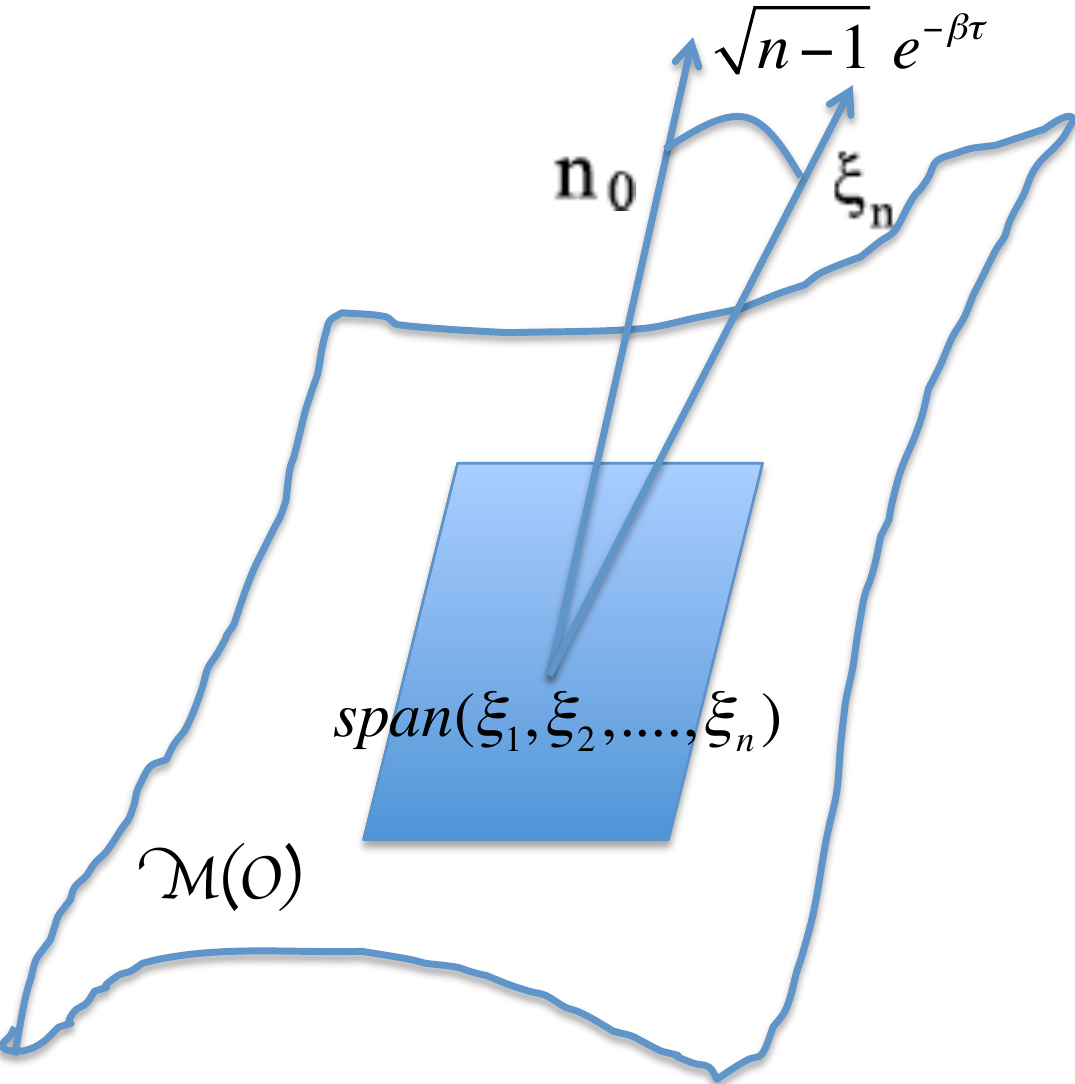}
  \caption{ The schematic showing decay in angle between the maximal expansion direction and normal direction for a repelling surface. }
 \label{fig_FTLE_rho_sche}
\end{figure}

The FTLE tries to estimate maximum local expansion along the trajectory, and can be used as a stability certificate for online stability monitoring. On the other hand, the repulsion ratio measures the rate of expansion along the normal direction. For normally repelling surfaces the direction of the maximum expansion and the normal direction align with each other asymptotically, which means that asymptotically the stability certificate obtained using both normal expansion rate and FTLE will be the same.  In order to obtain precise bounds of the alignment, we introduce the following notations - let $0 < \lambda_1 (x_0, \tau) \le \lambda_2 (x_0, \tau), \dots \lambda_n (x_0, \tau)$ be the eigenvalues of the matrix $\nabla \Phi (x_0, \tau) ^T \nabla \Phi (x_0, \tau)$, and $ \xi_1 (x_0, \tau), \xi_2 (x_0, \tau), \dots, \xi_n (x_0, \tau)$ be the corresponding eigenvectors. It can be noted that $\sigma^{\tau}(x_0) = \frac{1}{\tau} \log \lambda_n (x_0, \tau)$, and also $\xi_n (x_0, \tau)$ is the direction of the maximum expansion. The FTLE is the expansion along the  $\xi_n (x_0, \tau)$. In \cite{haller2011variational}, it was proved that  for a repelling hyperbolic surface the angle between the normal direction $\eta_0$ and the direction of maximum expansion $\xi_n (x_0, \tau)$ decays exponentially faster, where the rate of the decay is given by negative of the repulsion ratio of the repelling surface. Figure \ref{fig_FTLE_rho_sche} depicts the asymptotic alignment of the two vectors. It was shown in \cite{haller2011variational}, the sine of the angle between two vectors $\alpha_n$ satisfies the following condition, $\sin \alpha_n (\eta_0, \xi_n) \le \sqrt{n-1} e ^ {-\beta \tau} $. This alignment ensures the stability certificates FTLE and normal repulsion rate match with each other. \\ Apart from online stability monitoring, we were also interested  to compute finite time  stability boundaries. We can use normal expansion rate to get the finite time stability boundaries accurately. The stability boundaries thus identified from normally hyperbolic surface is called Lagrangian Coherent Structures (LCS) and are defined as following \cite{haller2011variational},
\begin{definition}
Assume that $\mathcal{M}(t)$ is a normally repelling material surface over $ t \in [ 0, \tau]$. We designate $\mathcal{M}(t)$ a repelling LCS	over	$[ 0, \tau]$	   if  its normal repulsion rate admits a point-wise non-degenerate maximum along  $\mathcal{M}(0)$ among all locally $C^1$-close material surfaces. We designate $\mathcal{M}(t)$ an attracting LCS over $[ 0, \tau]$ if it is a repelling LCS over $[ 0, \tau]$ in backward time. 
\end{definition}
The LCS are the local maximal ridges in formed in the scalar field of the normal repulsion rate. The stability boundaries can be identified by the LCS.  Whereas, FTLE can detect the stability boundaries only if certain additional conditions are satisfied. Under those conditions, FTLE ridges coincide with the LCS \cite{haller2011variational}. FTLE is a measure of maximum local expansion, and it can not differentiate between the normal and tangential expansion. But stability boundaries are essentially those surfaces, which are normally repelling. If a normally repelling surface separates two regions in state space either, and on either side of the stability boundary the trajectories have tangential expansion and normal contraction, FTLE ridges would  not be adequate to identify the stability boundaries. The reason behind this fact is that the FTLE field would show contraction for both the boundary and the neighboring region. On the other hand, the normal repulsion rate does not suffer from this fact as it would identify the fact the material surfaces on either side of the boundary are normally contracting. In \cite{haller2011variational} Theorem 13, the necessary and sufficient conditions for FTLE ridge to coincide with LCS and thus qualify as stability boundary is given. The conditions are - 1. $\lambda_{n-1} (x_0, \tau) \neq \lambda_n (x_0, \tau) > 1$, 2. $ \xi_n ( x_0, \tau) \perp T_{x_0} \mathcal{M} (0)$, 3. The matrix $L(x_0, \tau)$ ( which comprises of the $\lambda_i$'s and $\xi_i$'s )  is positive definite \cite{haller2011variational}.  We have used the normal expansion rate for both monitoring and stability boundary computation, whereas FTLE is used for online stability monitoring purposes. The FTLE is  computed using the model. 
\section{Model Free Algorithm to Compute FTLE}
Next, we introduce the model-free computation scheme from the time series to obtain approximate FTLE. In the following definition, we introduce the LE computation formula from the time series.
\begin{definition}
The Finite Time Lyapunov Exponent from time series for initial condition $x_0$ is defined as,
\begin{equation}  \label{LE_ts}
\lambda ( 0 , \tau ) : = \frac{1}{t} \log \frac{\parallel  \Phi (x_0,  \tau +\Delta t) - \Phi( x_0, \tau ) \parallel}{ \parallel  \Phi (x_0, \Delta t) - x_0 \parallel } , ~~ \tau >  \Delta t. 
\end{equation}
\end{definition}
The proposed method for LE computation is model free, and requires lesser computational effort when compared to normal expansion rate.   FTLE is positive if the trajectories diverge, and is negative if they converge. First, we propose a model-free scheme to compute LE from the time series. Then, we show that it is capable of identifying sources of local instabilities like presence of saddle points. This means if the trajectory comes close to a saddle point the LE starts to increase. Also, we demonstrate that FTLE has a cumulative effect, i.e. the local stability and instability contributions of various regions of the state space gets reflected as decrease or increase in the LE. This demonstrates its capability for being used as a stability certificate. 
\subsection{Model Free FTLE Based Stability Monitoring}
The Lyapunov Exponent is computed from the evolution of distance between actual and the delayed time series. Next, we present the following propositions, which relate LE with the transient stability problem. In Proposition \ref{LE_prop_1} we demonstrate that the LE would go negative for initial conditions inside the domain of attraction. This justifies the use of LE as a online stability certificate. 
\begin{proposition} \label{LE_prop_1}
For all initial point $x_0 \in A(x^s) \setminus x^s$, there exists a $T^*$ the finite time LE $ \lambda (x_0, 0 ,  t) < 0 $ for all $t \ge T^*$.
\end{proposition}

Proposition \ref{LE_prop_1} enables us to use LE as a stability certificate for online monitoring. The LE for the trajectories inside the domain of attraction will go negative, which can be used as a stability criterion. 

\section{Simulation Results}
In this paper, we develop theoretical foundation for  finite time transient stability monitoring of power systems. The theoretical foundation is based on tools from geometric theory of dynamical systems. In particular, we employ techniques from normally hyperbolic invariant manifold theory in the development of theoretical foundation.  We show that the normal expansion and contraction rate of co-dimension one manifold in the phase space can be used as a indicator for finite time transient stability. Furthermore extremum of these normal expansion and contraction scalar field can be used to identify stability boundary of the stable operating fixed point. Identification of stability boundaries or distance from the stability boundaries is used to determine margin of stability in transient. 
\subsection{Stability Boundary Computation}
\begin{figure}
    \centering
    \includegraphics[width=0.55\textwidth]{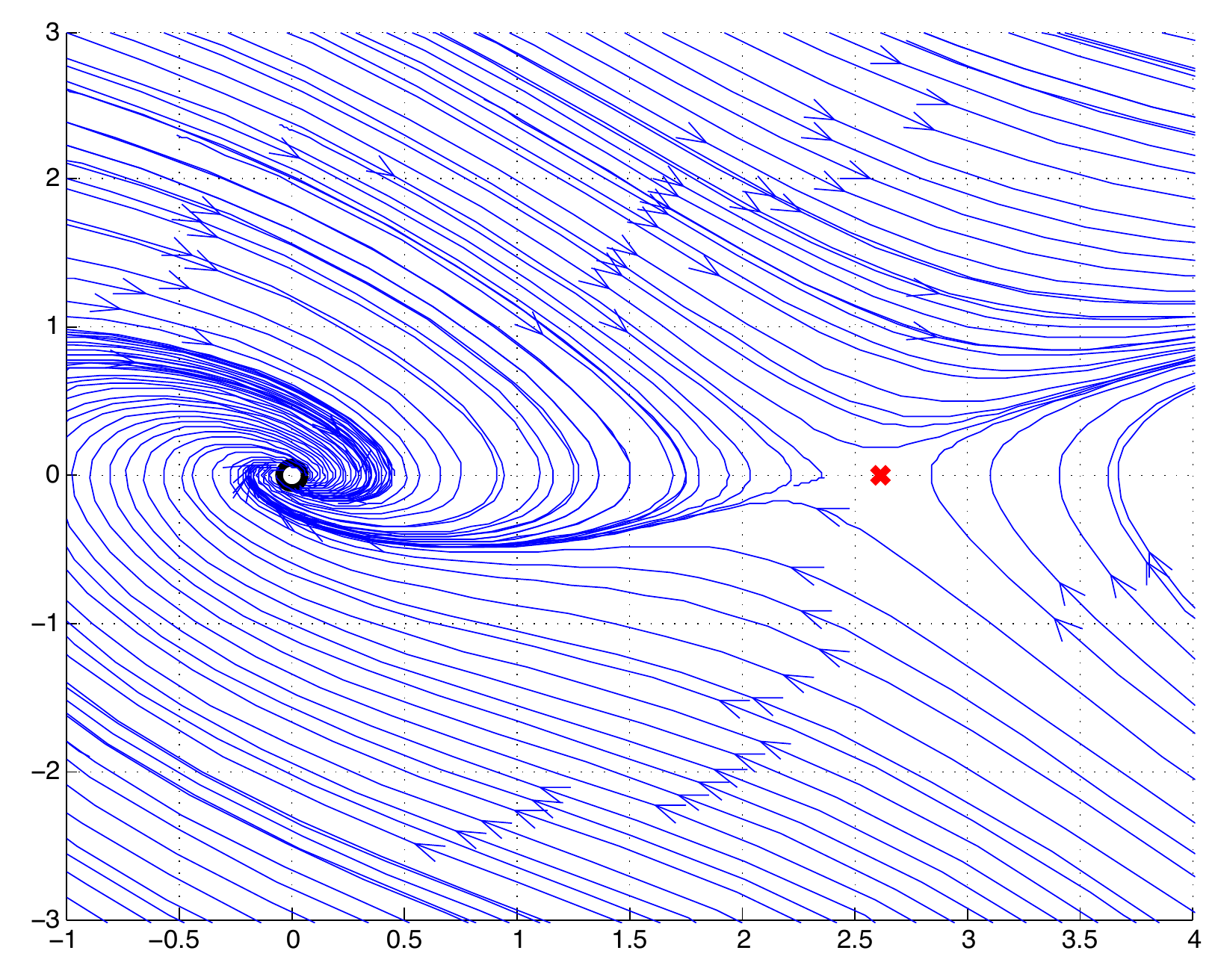}
  \caption{ Phase portrait of two generator system. }
 \label{fig_phase_portrait}
\end{figure}
\begin{figure}
    \centering
    \includegraphics[width=0.55\textwidth]{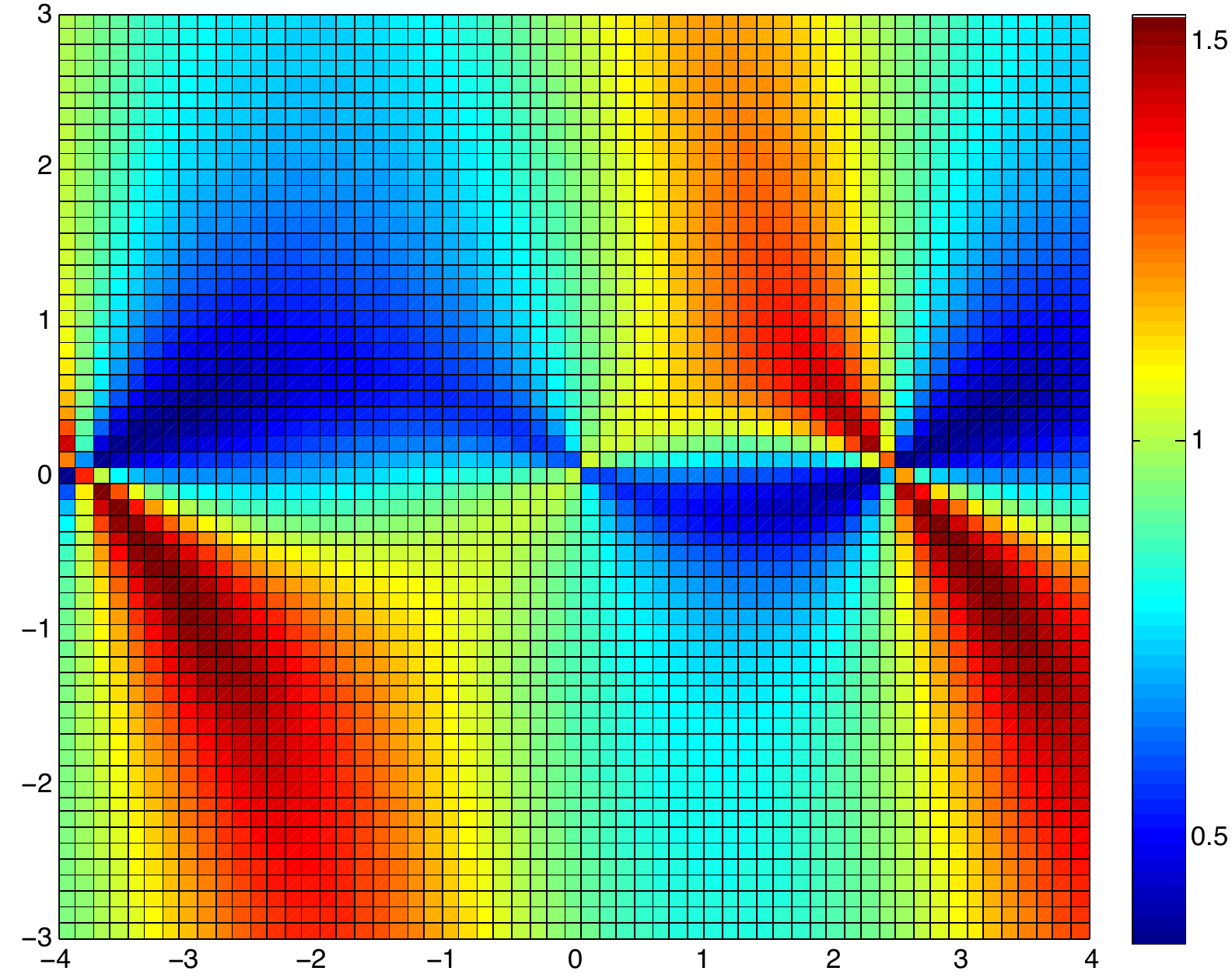}
  \caption{ Repulsion rate is plotted over a time window of $5$ secs. }
 \label{fig_LCS_5_sec}
\end{figure}
\begin{figure}
    \centering
    \includegraphics[width=0.55\textwidth]{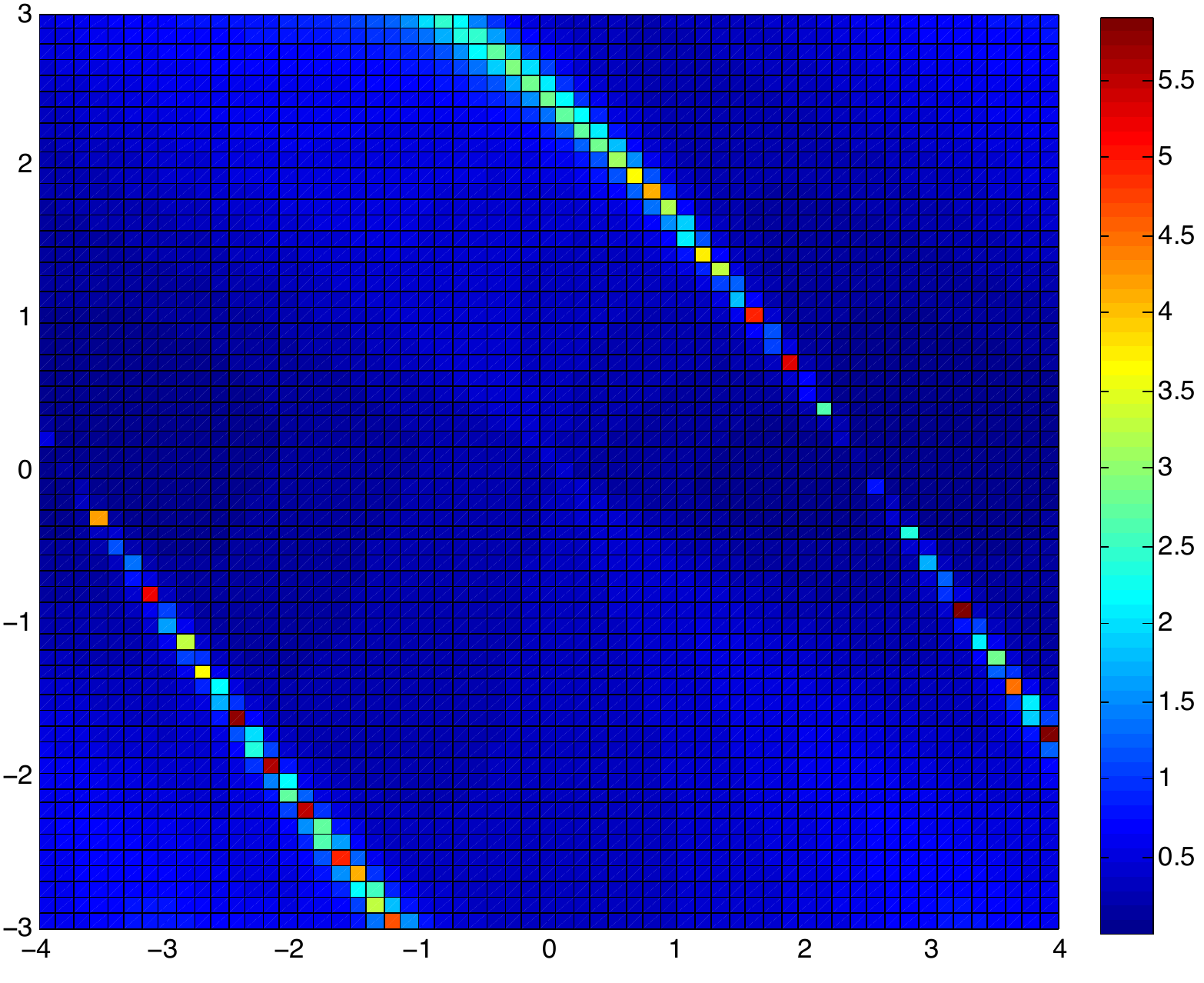}
  \caption{ Repulsion rate is plotted over a time window of $35$ secs. }
 \label{fig_LCS_35_sec}
\end{figure}
\begin{figure}
    \centering
    \includegraphics[width=0.55\textwidth]{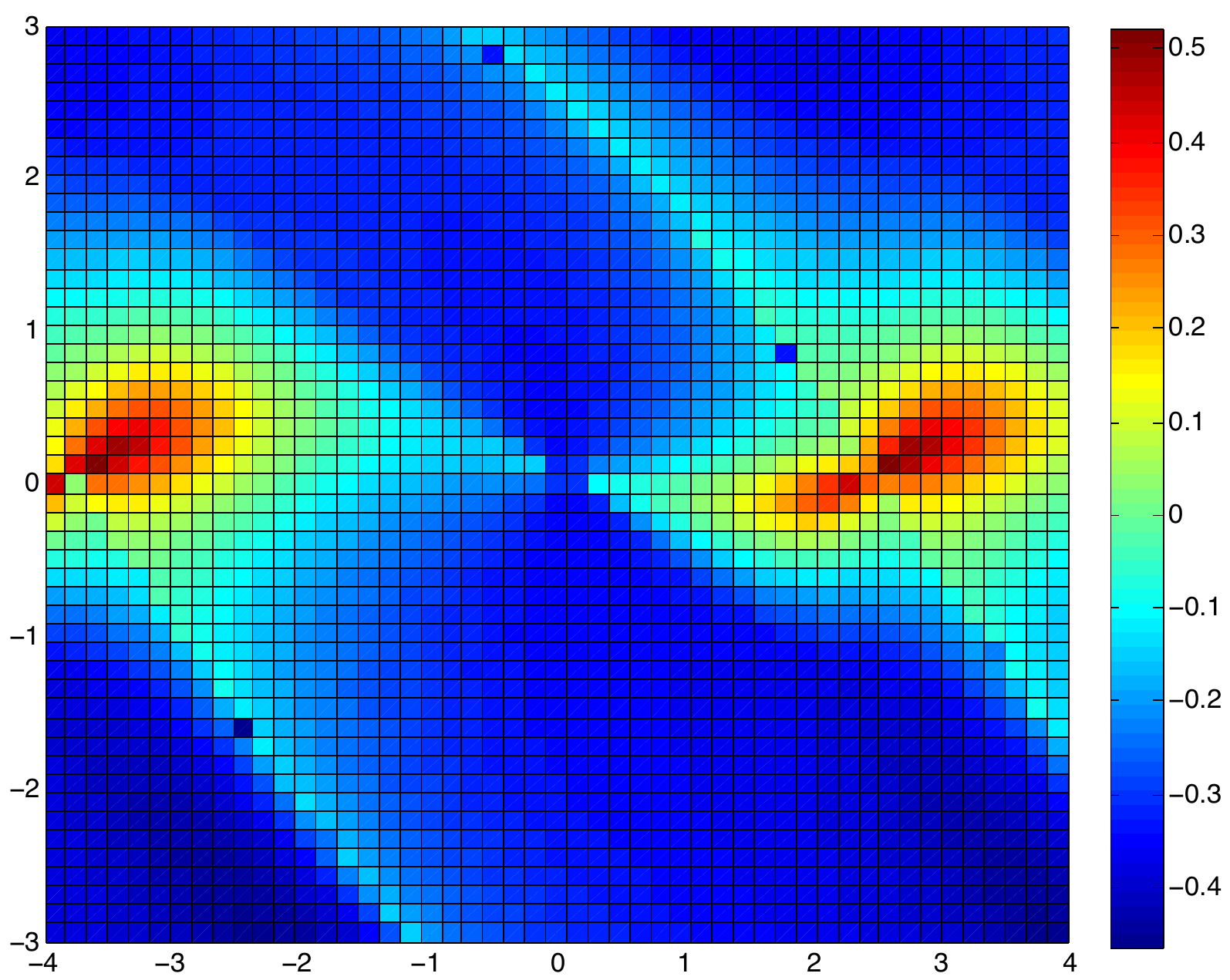}
  \caption{ LE is plotted over a time window of $5$ secs. }
 \label{fig_FTLE_5_sec}
\end{figure}
\begin{figure}
    \centering
    \includegraphics[width=0.55\textwidth]{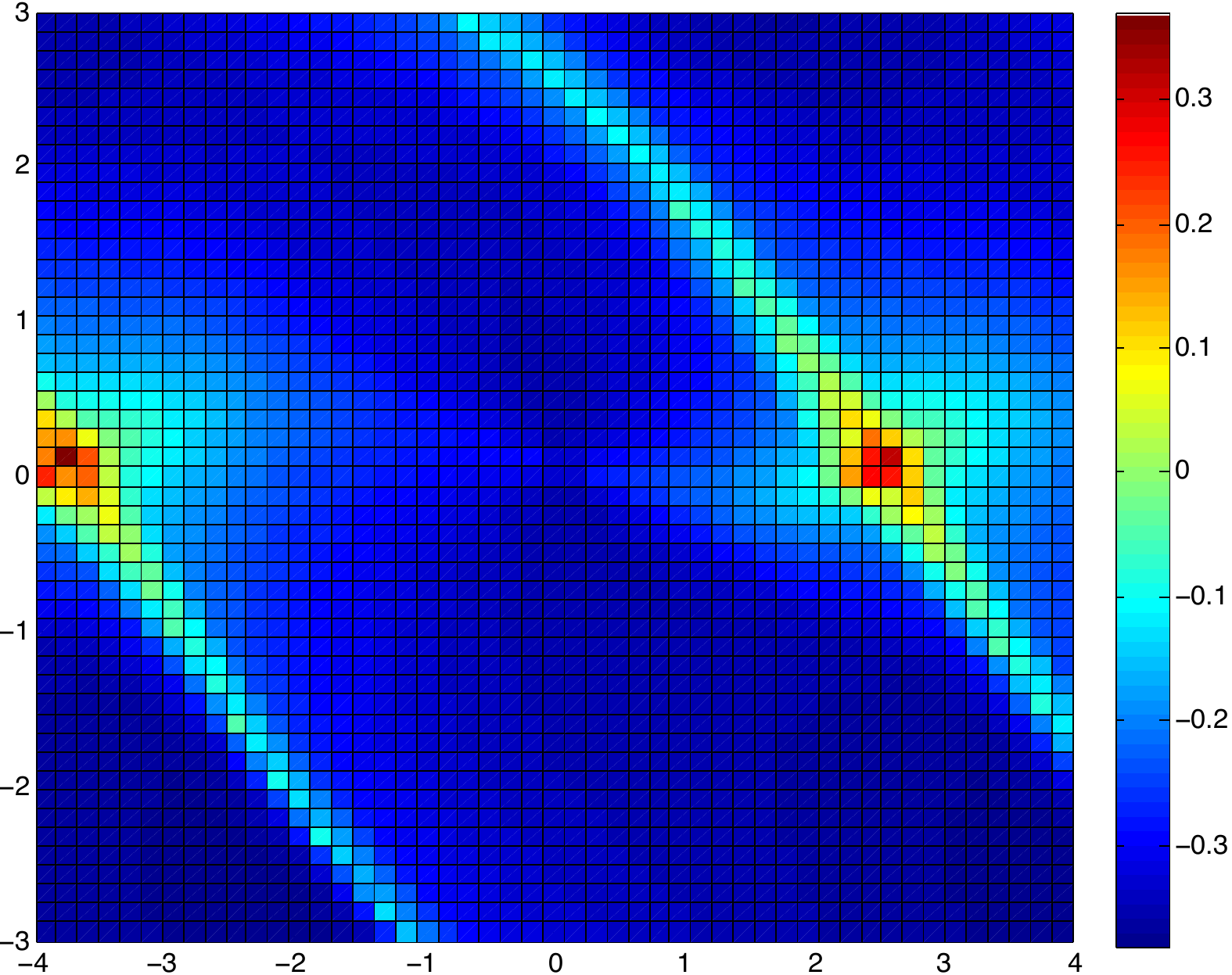}
  \caption{ LE is plotted over a time window of $35$ secs. }
 \label{fig_FTLE_35_sec}
\end{figure}
 For simulation purposes, we  consider the swing dynamics of $N+1$ generator system as the model, 
 \begin{align} \label{swing}
 \dot{\delta_i} &= \omega_i, \nonumber \\
 \dot{\omega_i} &= P_i - D_i {\omega}_i - \sum_{j=1}^{N+1} E_i E_j Y_{ij} sin \left ( \delta_i - \delta_j \right ), ~~~ i =1, \dots N.
 \end{align}
The $(N+1)^{th}$ generator is taken as the reference, by making $E_{N+1} = 1$, and $\delta_{N+1} = 0$. 
Let us define the state variable for the dynamical system as,
\begin{equation*}
 x : = [~\delta_1,~\dots,~ \delta_N, ~\omega_1, ~, \dots, ~\omega_N]^T \in \mathbb{R}^n. 
 \end{equation*}
This system corresponds to a specific set of values of the $P_i$'s, and also the system is at a stable equilibrium of the \eqref{swing}. Let, a fault occurs, at $t_f$, and is cleared at $ t = 0$. The source of the fault may be change in the electric power input at generator $i$, which is given by $P_i$, or short circuit, given by changes in the values $Y_{ij}$. Once the fault is cleared the system dynamics is changed according to the new parameters of the system. The phase portrait of two generator system is depicted in Fig. \ref{fig_phase_portrait}. The stable manifold of the type-$1$ saddle point forms the stability boundary. The stable equilibrium point is shown as a black circle and the type -$1$ saddle point is shown as a red cross. Figure \ref{fig_LCS_5_sec} shows the plot of the normal repulsion rate over a time window of $5$ sec, which is typical time interval for short term transient stability. Over this time window, it can be observed that a relatively larger region around the stable manifold of the type-$1$ saddle point show normal repulsion more than $1$. The region, which has more than $1$ repulsion rate, indicates the finite time stable region. This can not be captured by the asymptotic approaches. It can be observed from Fig. \ref{fig_LCS_35_sec}, as the time window is extended to $35$ sec, the stable manifolds of the type-$1$ saddle point emerges as the material surface with positive repulsion rate, which matches with the asymptotic analysis. It demonstrates, our approach gives a concrete way to provide the region of stability over a finite time window. Figures \ref{fig_FTLE_5_sec}, and \ref{fig_FTLE_35_sec} show the plot of FTLE over  time intervals of $5$, and $35$ sec. It can be observed that the FTLE can locally detect the saddle point. FTLE can not detect globally the stability boundaries, as it becomes negative along the stable manifold of the saddle point. Whereas, near the saddle point, the unstable manifold has locally positive FTLE. This is because the unstable manifold of the saddle is locally tangentially expanding. Also, we have observed the ridges in the normal repulsion rate or FTLE field preserve themselves under small perturbations or noise. Next, we demonstrate with simulations on $39$ bus system, that the normal expansion rate and LE can be used for online stability monitoring.
\subsection{Online Transient Stability Monitoring  for New England $39$ Bus System}
\begin{figure}
    \centering
    \includegraphics[width=0.55\textwidth]{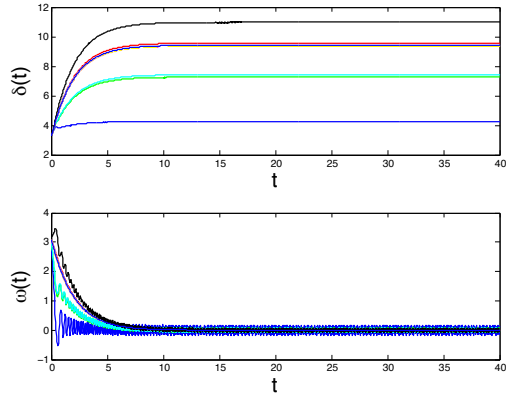}
  \caption{ Angle and frequency time series for stable dynamics in $39$ bus system. }
 \label{fig_39_39_bus_ts_stable}
\end{figure}
\begin{figure}
    \centering
    \includegraphics[width=0.55\textwidth]{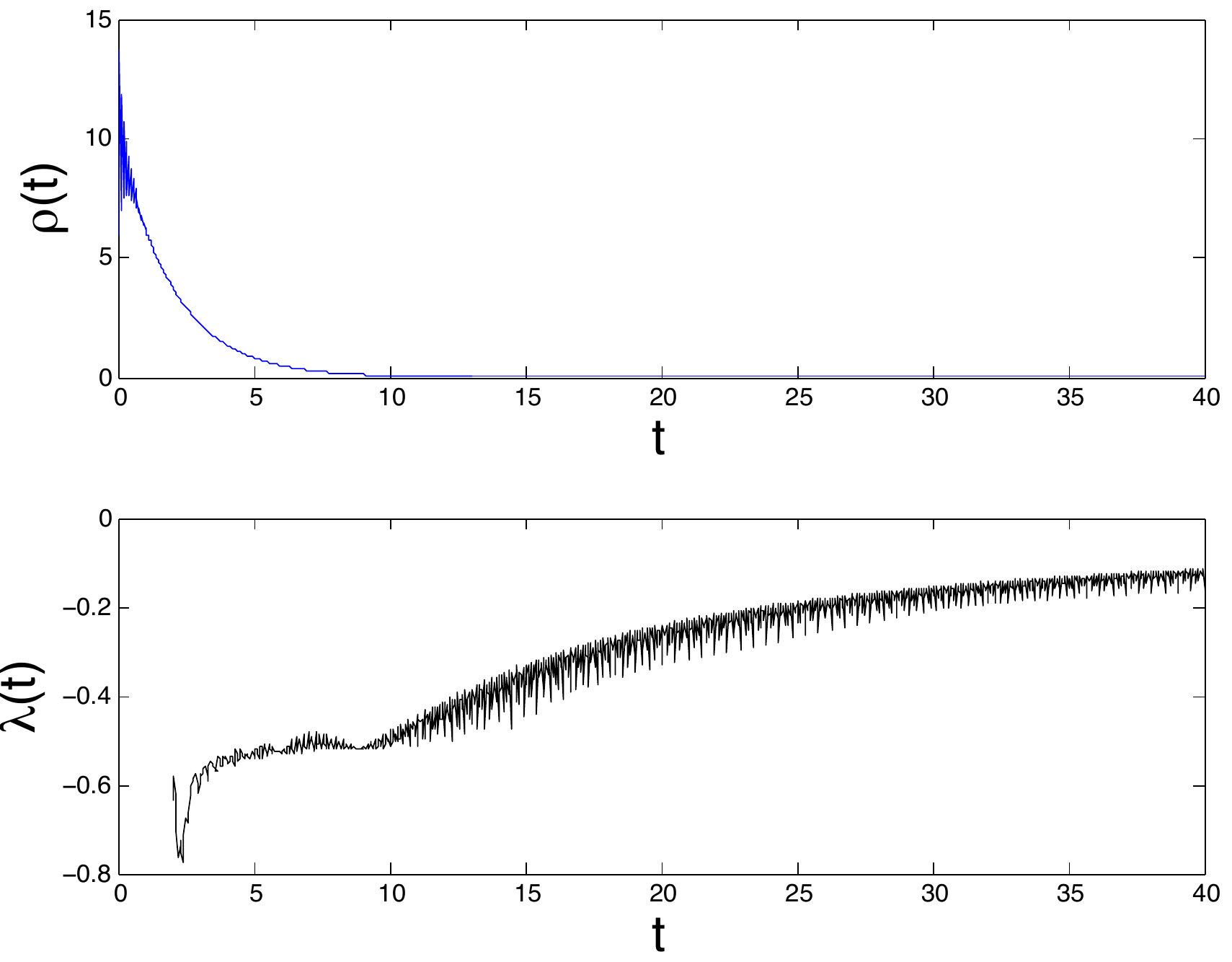}
  \caption{ Evolution of stability certificates $\rho(t)$, and $\lambda(t)$ for stable dynamics in $39$ bus system. }
 \label{fig_39_39_bus_rho_stable}
\end{figure}
\begin{figure}
    \centering
    \includegraphics[width=0.55\textwidth]{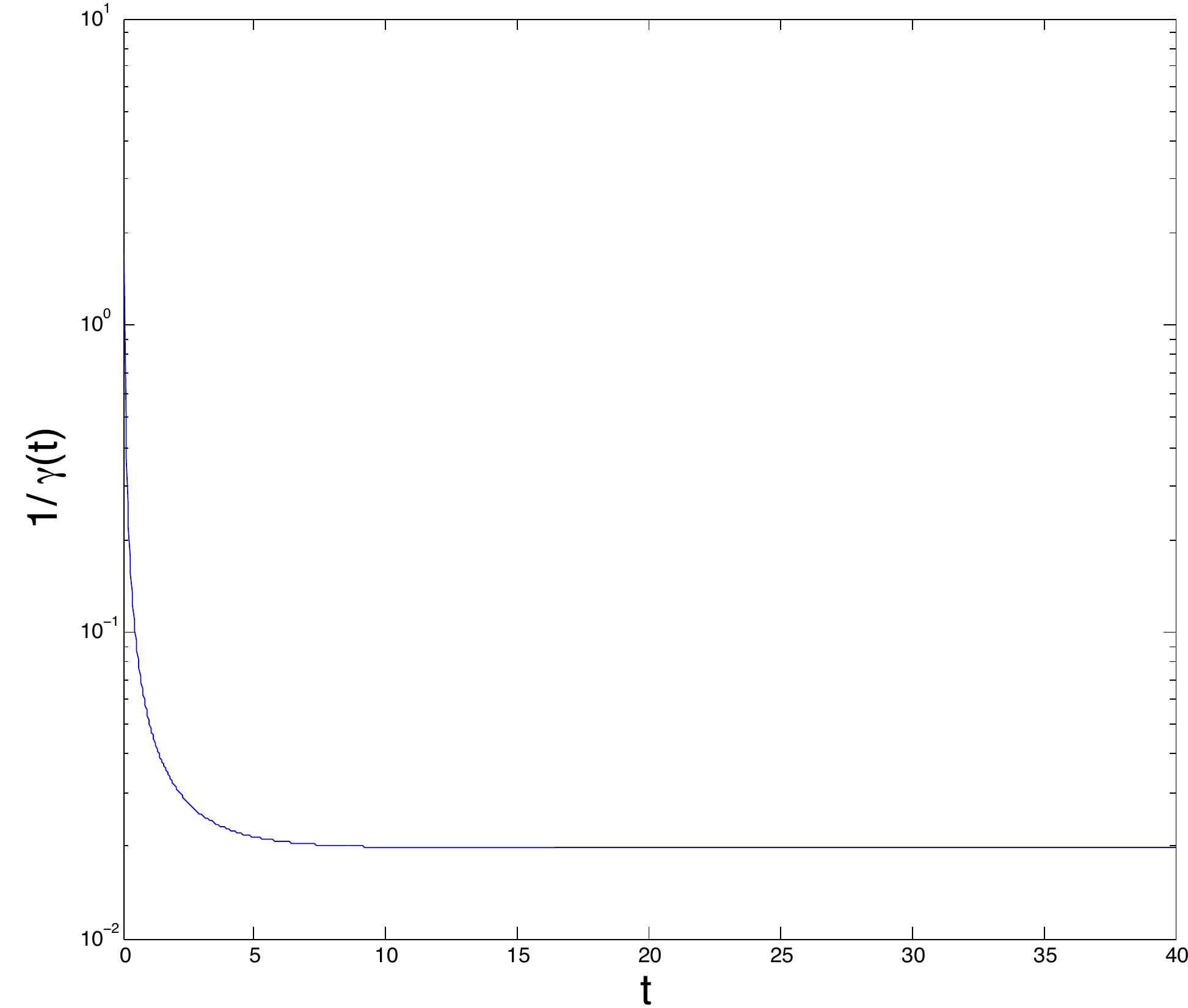}
  \caption{ Evolution of stability margin $\frac{1}{\gamma(t)}$ for stable dynamics in $39$ bus system. }
 \label{fig_39_39_bus_gamma_stable}
\end{figure}
\begin{figure}
    \centering
    \includegraphics[width=0.55\textwidth]{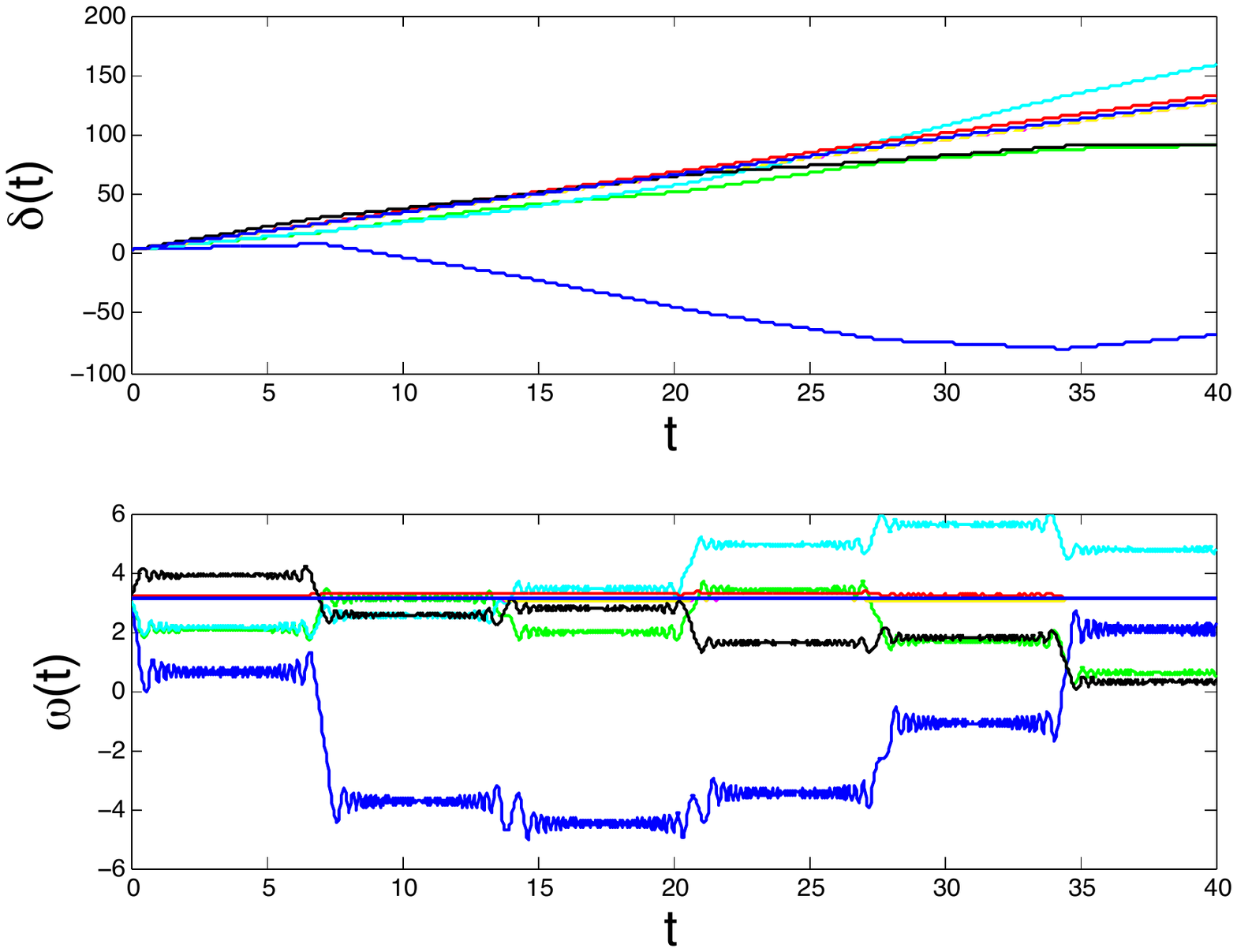}
  \caption{ Angle and frequency time series for unstable dynamics in $39$ bus system. }
 \label{fig_39_39_bus_ts_unstable}
\end{figure}
\begin{figure}
    \centering
    \includegraphics[width=0.55\textwidth]{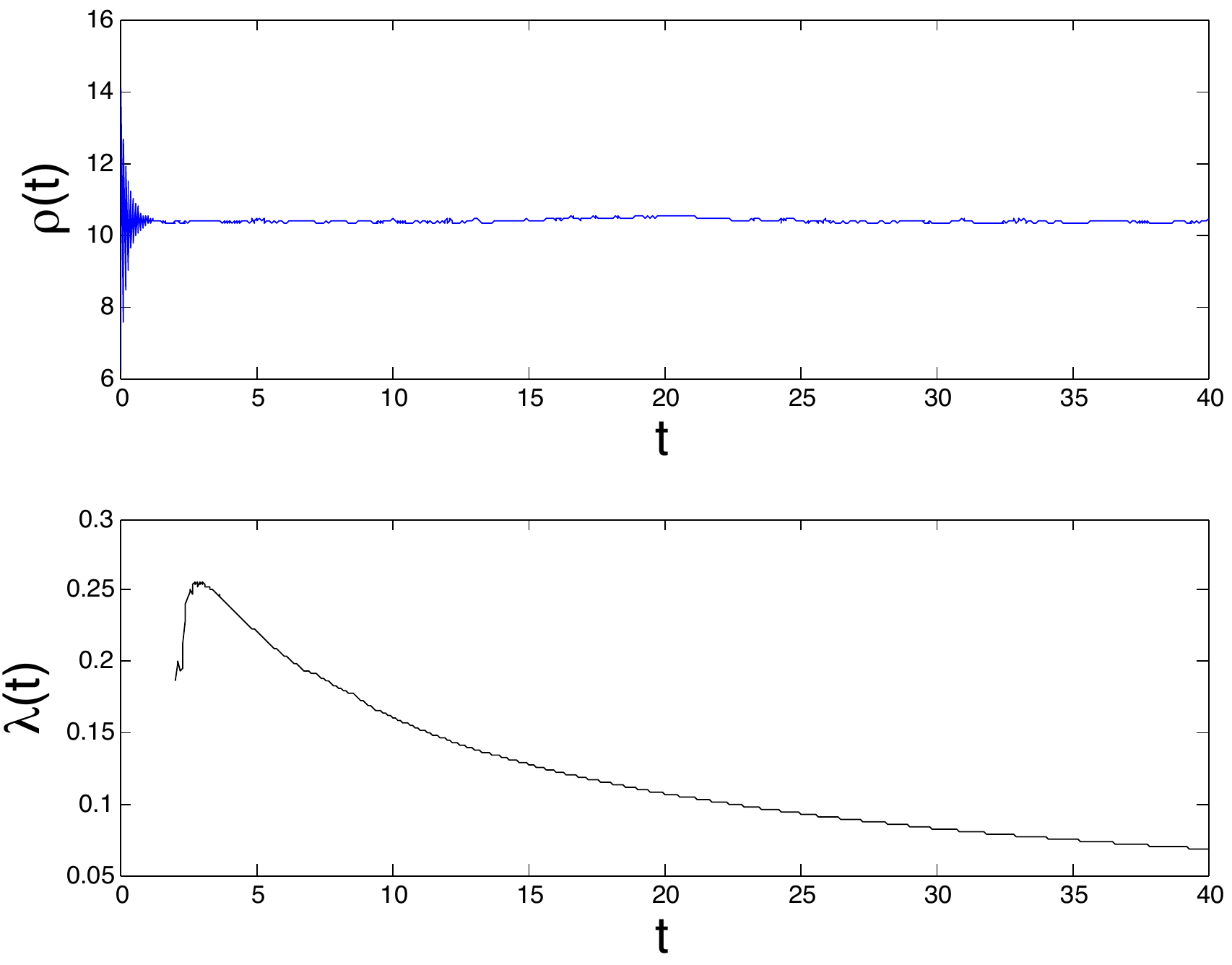}
  \caption{ Evolution of stability certificates $\rho(t)$, and $\lambda(t)$ for unstable dynamics in $39$ bus system. }
 \label{fig_39_39_bus_rho_unstable}
\end{figure}
 \begin{figure}
    \centering
    \includegraphics[width=0.55\textwidth]{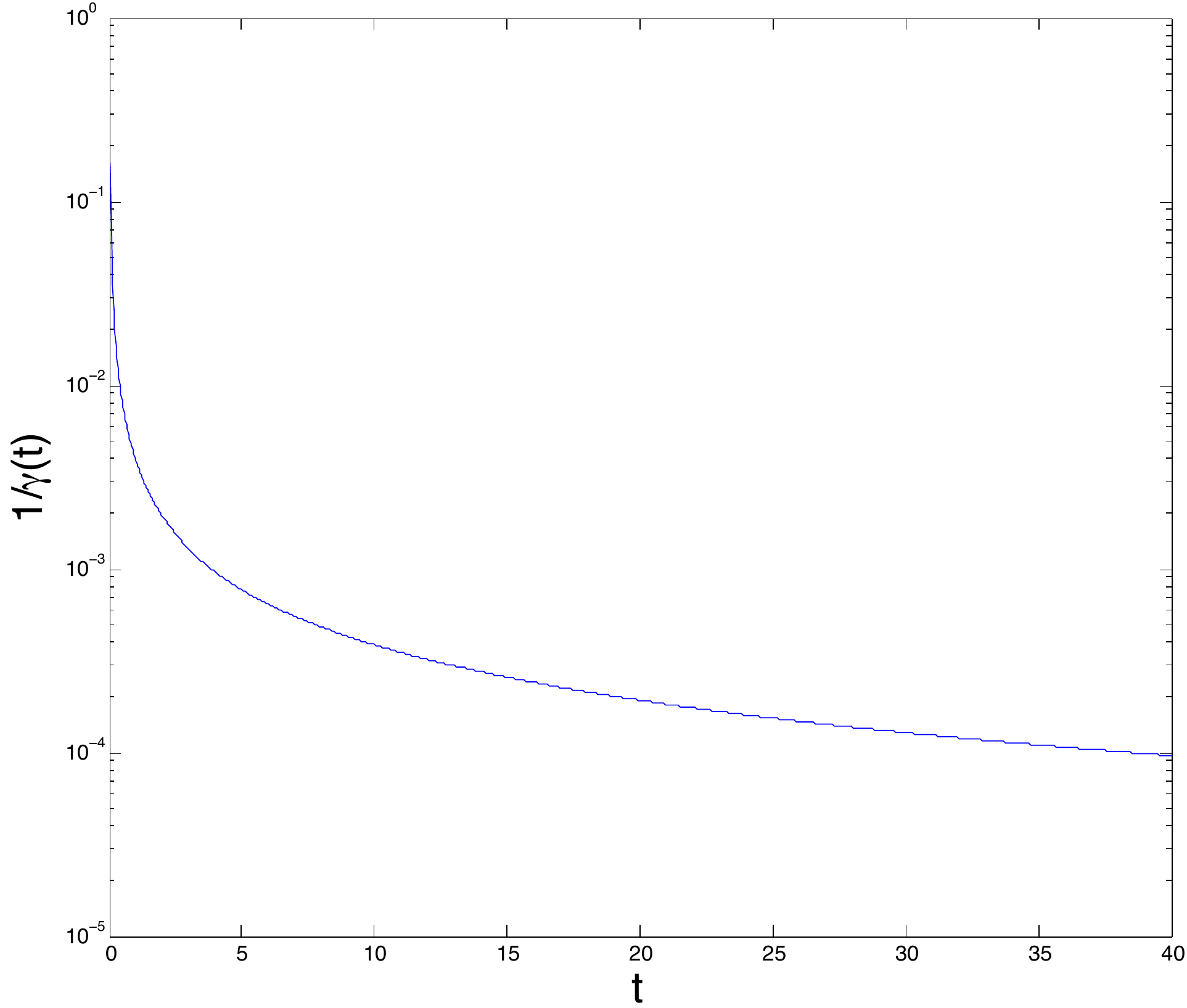}
  \caption{ Evolution of stability margin $\frac{1}{\gamma(t)}$ for unstable dynamics in $39$ bus system. }
 \label{fig_39_39_bus_gamma_unstable}
\end{figure}
In this subsection, we show application of normal expansion rate  $\rho(t)$  , and LE  $\lambda (t)$   as stability certificates in online monitoring. Also, we compute stability margins based of normal expansion rate, and LE, which are denoted as $\frac{1}{\gamma (t)}$, and $\frac{1}{\theta (t)}$. The stability certificates, as well as margins can further be used to generate alarms, and to generate appropriate control actions. The computational time for LE is less than that of the normal expansion rate. On the other hand, the normal expansion rates can more accurately detect the global stability boundaries. Thus, depending on the available computational infrastructure, and desired accuracy, the right stability certificate could be chosen. \\ $~~~$ We present the simulation results for the swing dynamics of the New England $39$ bus system. New England $39$ bus system has $10$ generators, and it is a reduced model for the power grid of New England and part of Canada. We have used the normal rate of expansion $\rho(t)$, $\lambda$ as a stability certificates.  The normal expansion rate is computed according to \eqref{rho_formula}, and LE is computed as described in \eqref{LE_ts}. The $\Delta t$ for LE computation is chosen as $1.5$ sec.  \\$ ~~~$ The governing equations of the rotor angles swing dynamics are  as following \cite{kundur2001power_book}, 
\begin{align*}
\dot \delta_i (t) & = \omega_i (t), \\
\frac{H_i}{\pi f_s} \omega_i &  = - D_i \omega_i + P_{mi} - G_{ii} E_{i}^2 \\ 
& -  \sum_{j=1, i \neq j}^{10} E_i E_j \left ( G_{ij} \cos(\delta_i - \delta_j ) + B_{ij} \sin(\delta_i - \delta_j )  \right ),
\end{align*}
where, $\delta_i$, and $\omega_i$ are the angle and the frequency of the $i^{th}$ generator. The values of the system parameter can be found in \cite{mat_power}. We create a stable and unstable scenario by tuning the damping parameter. Figure \ref{fig_39_39_bus_ts_stable} shows the system dynamics for the stable scenario (damping $ D_i = 0.5 $). It can be observed from Fig. \ref{fig_39_39_bus_rho_stable} that the normal expansion rate $\rho(t)$ decays exponentially faster and converges to $0$, showing stable behavior, and the LE stays below $0$, prescribing stability. Figure \ref{fig_39_39_bus_gamma_stable} shows the evolution of the stability margin $\frac{1}{\gamma(t)}$. It can be observed that, in this case the stability margin converges to a constant value as the system goes stable.  \\ On the other hand, Figure \ref{fig_39_39_bus_ts_unstable} shows the angles, and frequencies, when the system is unstable (damping $ D_i = 0.01 $). It can be observed from Fig. \ref{fig_39_39_bus_rho_unstable} that the normal expansion rate $\rho(t)$ stays above $1$, certifying the instability. Also, the LE stays above $0$, showing instability. Figure \ref{fig_39_39_bus_gamma_unstable} shows the evolution of the stability margin $\frac{1}{\gamma(t)}$ for the unstable. It can be observed that, in this case the stability margins show a monotone decreasing trend, which indicates unstable behavior.  The simulations testifies that the normal expansion rate and LE can be used as a online stability monitoring tool.
\section{Conclusion}
In this paper, we propose a finite time stability analysis tool based on the Theory of normally hyperbolic surfaces. We have related the stability boundaries of the transient dynamics to the normally hyperbolic repelling surfaces. Our proposed method can prescribe stability region based on the time window of interest, which is very useful for finite time transient stability problem. The normal repulsion rate can also be used for online stability monitoring. We have also proposed a LE based model-free stability monitoring scheme for fast real time applications.
\section{Appendix}
Next, we outline the proof of Theorem \ref{norm_hyp_1}. Theorem \ref{norm_hyp_1} aims at showing that there would be a region around the stable manifold, which would be normally repelling over a finite window of time. Before outlining the details of the proof, we would try to put forth the intuition behind the successive constructions. The stable manifold forms a co-dimensional manifold, which separates the trajectories inside and outside the domain of attraction to its either side. If we consider two points on the opposite sides of the stable manifold, separated by a small distance. The line joining these two points would be normal to the stable manifold. Trajectories emanating from these two points would evolve in different fashion - the point inside the domain of attraction would eventually converge to the stable fixed point and the other one would stay outside the domain of attraction. As a result, the distance between these two points would be beyond a threshold. If the initial distance between the two points have taken very small, the ratio between the initial and final distance could be made arbitrarily large. We would use this fact to show that the stable manifold would emerge as a normally repelling surface over a finite time window. Next, we would outline the technical details, relating to the proof. \\ 
\begin{proof}
Let us consider the stable manifold of the type-$1$ saddle point $x^e$, which is denoted by $W_s(x^e)$. We need to show for every $T>0$ and $x \in W_s (x^e) $, there exists an $\epsilon > 0$ such that every point in the set $\hat W_u (x) $ has repelling hyperbolic material surface, where,  
\begin{align*}   
\hat W_u (x)  & := \{ ~ \hat x  ~ | ~  \parallel x - \hat x \parallel < \epsilon , ~   x \in W_s (x^e)  \}. 
\end{align*}
We can construct two points  $x_0 - \epsilon n_0 (x_0) $, and $x_0 + \epsilon n_0(x_0)$, where $n_0(x_0)$ is an arbitrary unit vector - these would be the two points on either side of the stable manifold.  Next, we consider the evolution of distance between  $\Phi \left ( x_0 - \epsilon n_0 (x_0), t \right ) $, and $\Phi \left ( x_0 + \epsilon n_0 (x_0), t \right)$, which would finally lead us to the proposition. Now let us define,
\begin{align*} 
e(t, \epsilon) & := \Phi(x_0 + \epsilon n_0 (x_0), t ) - \Phi(x_0 - \epsilon n_0 (x_0), t ) \\
& = 2 \epsilon \frac{\partial \Phi}{\partial x } (x_0 - \epsilon n_0 (x_0), t ) n_0 (x_0) + \mathcal{O}(\epsilon^2).
\end{align*}
where, $\mathcal{O}(\epsilon^2)$ can be ignored for sufficiently small $\epsilon$. Now, for a given $T > 0$, 
\begin{align} 
e(T, \epsilon) & =  2 \epsilon \frac{\partial \Phi}{\partial x } (x_0 - \epsilon n_0 (x_0), T ) n_0 (x_0) \nonumber \\
e' (T , \epsilon ) e(T , \epsilon ) & = 4  \epsilon^2 n'_0 (x_0) \mathcal M (x_0, n_0, \epsilon, T) n_0 (x_0), \label{sq_norm}
\end{align}
where, \[ \mathcal M (x_0, n_0, \epsilon, T) := \left (  \frac{\partial \Phi}{\partial x } (x_0 - \epsilon n_0, T )  \right )' \frac{\partial \Phi}{\partial x } (x_0 - \epsilon n_0, T ) \] is a positive definite matrix. The singular value decomposition of the matrix $\mathcal M (x_0, n_0, \epsilon, T)$ would give us,
\begin{equation*}
\mathcal M (x_0, n_0, \epsilon, T) = U' (x_0, n_0, \epsilon, T) ~\Sigma (x_0, n_0, \epsilon, T) ~ U (x_0, n_0, \epsilon, T), \label{svd}
\end{equation*}
where, $U (x_0, n_0, \epsilon, T)$ is a unitary matrix with $i^{th}$ column $u_i (x_0, n_0, \epsilon, T)$, and $\Sigma (x_0, n_0, \epsilon, T) = diag\left (\Lambda_1 (x_0, n_0, \epsilon, T) , \dots, \Lambda_n (x_0, n_0, \epsilon, T)  \right)$, where $\Lambda_1 \ge \Lambda_2 \ge \dots \ge \Lambda_n > 0$. This gives us,
 \begin{align}
&  \mathcal M (x_0, n_0, \epsilon, T) = \nonumber \\
&\sum_{i=1}^N \Lambda_i (x_0, n_0, \epsilon, T) u_i (x_0, n_0, \epsilon, T)  u_i'  (x_0, n_0, \epsilon, T) . \label{svd_2}
 \end{align}
 \begin{figure}
    \centering
    \includegraphics[width=0.45\textwidth]{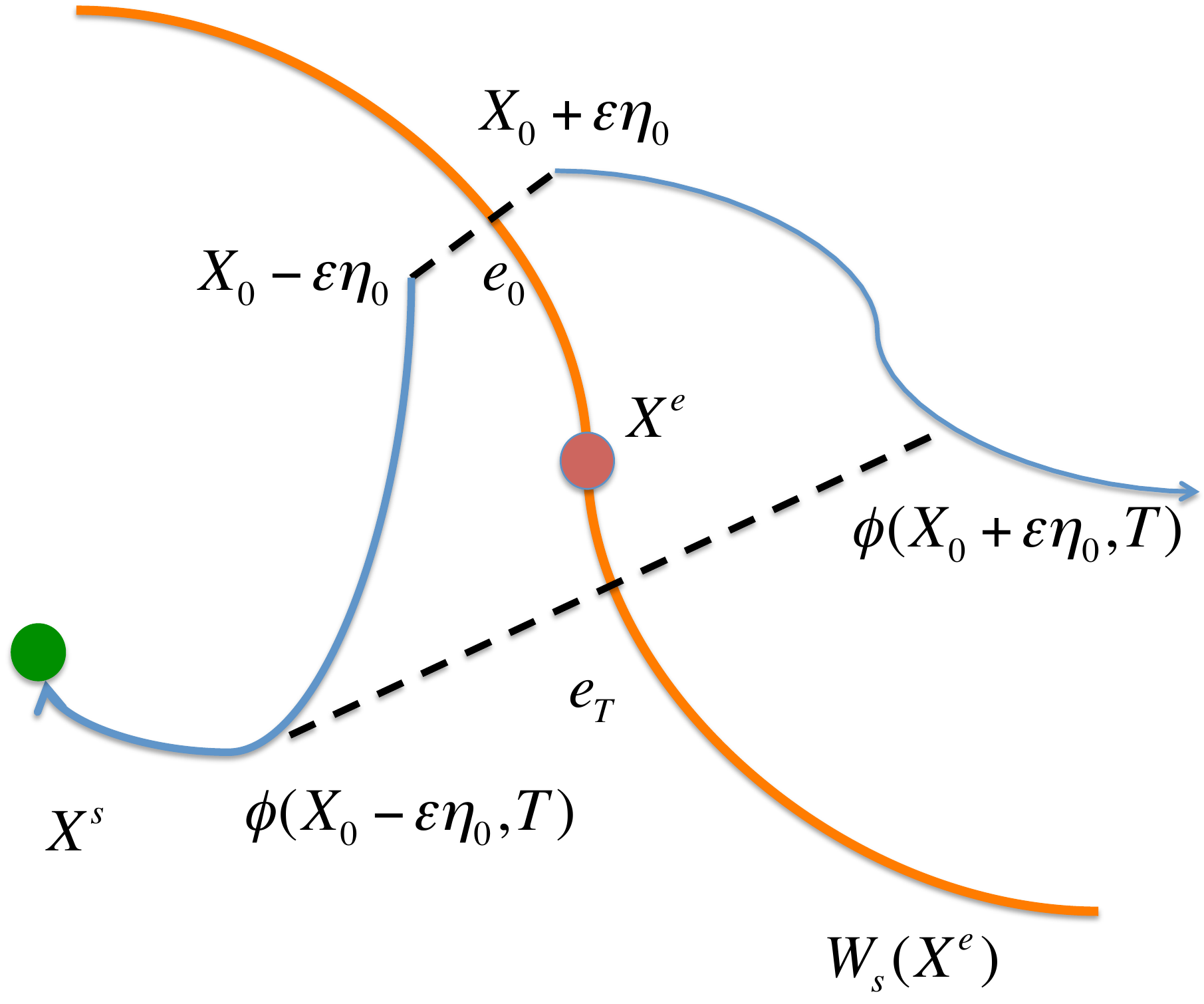}
  \caption{ Schematic of proof of Theorem \ref{norm_hyp_1}. }
 \label{fig_norm_hyp_1}
\end{figure} 
 Combining \ref{sq_norm} and \ref{svd_2}, 
 \begin{align*}
 & e' (T , \epsilon) e(T , \epsilon)  = \\
 &4  \epsilon^2  \sum_{i=1}^N \Lambda_i (x_0, n_0, \epsilon, T) n'_0 u_i (x_0, n_0, \epsilon, T)  u_i'  (x_0, n_0, \epsilon, T) n_0. 
 \end{align*} 
By appropriately selecting  $\epsilon$, 
\[ e' (T , \epsilon) e(T , \epsilon)  > \mathcal {K}(T,\epsilon ) e' (0 , \epsilon) e(0 , \epsilon), \]
Next, we would inspect the quantity $\mathcal{K} (T, \epsilon)$. It can be noted that, $e(0 , \epsilon) = 2 \epsilon n_0, \text {and} ~ e'(0 , \epsilon) e(0 , \epsilon) = 4 \epsilon^2 n'_0 n_0 =  4 \epsilon^2 $, as $ n_0' n_0 = 1$. Now, the quantity, $e(T ,\epsilon) $ would be more than a finite $\hat K$ after a finite $T$, as  $e(T , \epsilon)$ shows the distance between two points on the either side of the separatrix. Hence $ \frac{e' (T , \epsilon) e(T , \epsilon)}{ 4 \epsilon^2}> \mathcal{K} (T, \epsilon)$. So, by choosing sufficiently small $\epsilon$ we can have an arbitrary large $\mathcal{K} (T, \epsilon)$.
We do the further simplification of the SVD,  
\begin{align*}
 & 4  \epsilon^2  \sum_{i=1}^n \Lambda_i (x_0, n_0, \epsilon, T) n'_0 u_i (x_0, n_0, \epsilon, T)  u_i'  (x_0, n_0, \epsilon, T) n_0\\ 
 & > 4  \epsilon^2 \mathcal{K} (T ,\epsilon), \\
&  \sum_{i=1}^n \Lambda_i (x_0, n_0, \epsilon, T) n'_0 u_i (x_0, n_0, \epsilon, T)  u_i'  (x_0, n_0, \epsilon, T) n_0 \\
& > \mathcal{K} (T ,\epsilon)  , \\
& \sum_{i=1}^n \Lambda_i (x_0, n_0, \epsilon, T) \left ( n'_0 u_i (x_0, n_0, \epsilon, T) \right )^2 \\
 & > \mathcal{K} (T ,\epsilon)  . 
\end{align*}
By choosing, $n(x_0) = u_n$, we get, $\Lambda_n (x_0, n_0, \epsilon, T)  > \mathcal{K} (T ,\epsilon) $. 
For the unit normal vector $\eta_0 (x_0)$, we have, 
\begin{align*} 
& \frac{1}{  \sum_{i=1}^n \frac{1}{\Lambda_i (x_0, n_0, \epsilon, T)} \left ( \eta'_0 u_i (x_0, n_0, \epsilon, T) \right )^{ 2} }  = \frac{1}{  \sum_{i=1}^n \frac{1}{\Lambda_i (x_0, n_0, \epsilon, T)} \eta'_0 \eta_0 }\\
& =  \frac{1}{  \sum_{i=1}^N \frac{1}{\Lambda_i (x_0, n_0, \epsilon, T)} }  > \frac{\Lambda_n}{n} >  \frac{\mathcal{K} (T ,\epsilon, T)}{n} .   
\end{align*}  
The $\epsilon$ can be made sufficiently small such that, $  \mathcal{K} (T ,\epsilon) > n $. This will give, 
\begin{align}
\frac{1}{  \sum_{i=1}^n \frac{1}{\Lambda_i (x_0, n_0, \epsilon, T)} \left ( \eta'_0 u_i (x_0, n_0, \epsilon, T) \right )^{ 2} } & > 1, \\
 \frac{1}{  \eta_0^T  \left ( \left( \nabla  \Phi(x_0 - \epsilon n_0 , T)) \right)^T \nabla  \Phi(x_0 - \epsilon n_0 , T)  \right ) ^{-1} \eta_0 } & > 1, \\
 \rho(x_0 - \epsilon n_0, T) > 1.
\end{align}
Similarly, to prove $ \rho(x_0 +  \epsilon n_0, T) > 1$, we use,
\[ e(T, \epsilon)  =  - 2 \epsilon \frac{\partial \Phi}{\partial x } (x_0 +  \epsilon n_0 (x_0), T ) n_0 (x_0)  \] and use the successive steps as it is. \\ 
Hence, the proof.
\end{proof}
 \begin{figure}
    \centering
    \includegraphics[width=0.45\textwidth]{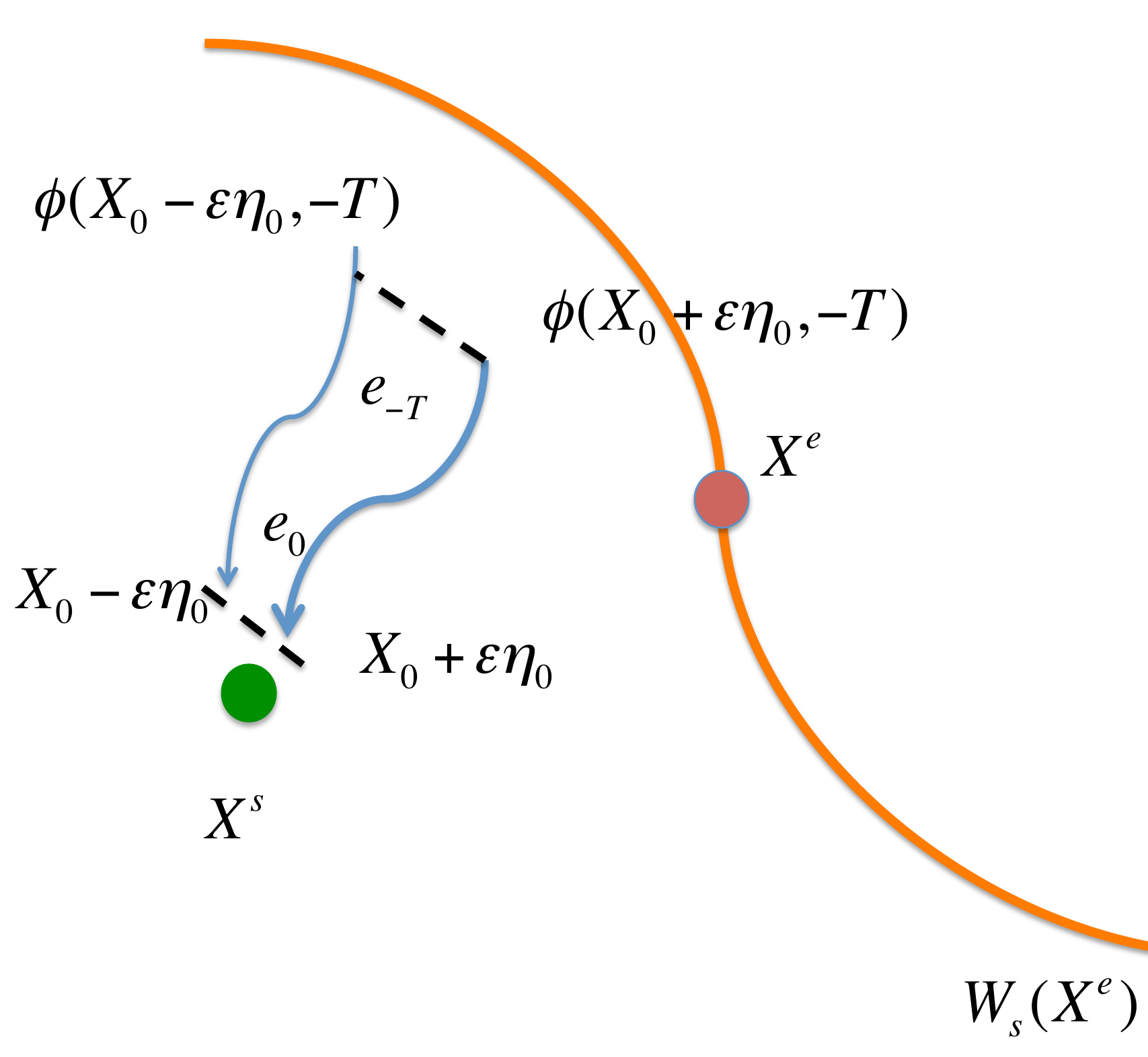}
  \caption{ Schematic of proof of Theorem \ref{norm_hyp_2}. }
 \label{fig_norm_hyp_2}
\end{figure}
We outline the proof of Theorem \ref{norm_hyp_2}. Here also the proof is along the same lines as the last one. We take two points within the domain of attraction, and use the fact that the they would eventually converge to the stable fixed point. Figure \ref{fig_norm_hyp_2} captures the essence of the construction. \\
\begin{proof}
Let us first consider the stable manifold of the type-$1$ saddle point $x^e$, which is denoted by $W_s (x^e) $. We need to show for every $T>0$ and $x \in W_s (x^e)  $, there exists an $\epsilon > 0$ such that every point in the set $\hat W_s (x) $ has attracting hyperbolic material surface, where,  
\begin{align*}   
\hat W_s (x)  & := \{ ~ \hat x  ~  | \parallel x - \hat x \parallel > \epsilon ,   ~ x \in W_s   (x^e) \}. 
\end{align*}
The trajectories in domain of attraction $A(x^s)$ would asymptotically converge to the stable fixed point $x^e$. Stability of the fixed point $x^s$ gives for every $T$ and $\delta_1 > 0$ there exists a $\delta_2 > 0$ such that, for any point $x$ in the $\delta_2$ neighborhood of $x^s$,  $\parallel \Phi(x,T) - x^s \parallel < \delta_1 $.  The attracting material surfaces are defined in terms of the reverse flow. The stability condition implies, for every $T$ and $\delta_1 > 0$ there exists a $\delta_2 > 0$ such that, for any point $z$ such that $\parallel \Phi(x,-T) - x^s \parallel$ ,  satisfies $\parallel x - x^s \parallel < \delta_1$. Now, let us choose two points $x_0$, and $x_0+\epsilon_0 n_0(x)$ such that $\epsilon_0 < \frac{\delta_1}{2}$, and $n_0(x)$ is an arbitrary vector.  
Next, we consider the evolution of distance between  $\Phi \left ( x_0 - \epsilon_0 n_0 (x_0), -t \right) $, and $\Phi \left ( x_0 + \epsilon_0 n_0 (x_0), -t \right )$. Now let us define,
\begin{align*} 
e( - t, \epsilon_0) & := \Phi(x_0 + \epsilon_0 n_0 (x_0), - t ) - \Phi(x_0 - \epsilon_0  n_0 (x_0),  - t ) \\
& = 2 \epsilon_0 \frac{\partial \Phi}{\partial x } (x_0 - \epsilon_0 n_0 (x_0), - t ) n_0 (x_0) + \mathcal{O}(\epsilon_0^2).
\end{align*}
where, $\mathcal{O}(\epsilon_0^2)$ can be ignored for sufficiently small $\epsilon_0$. Now, for a given $T > 0$, 
\begin{align} 
e(- T, \epsilon_0) & =  2 \epsilon_0  \frac{\partial \Phi}{\partial x } (x_0 - \epsilon_0  n_0 (x_0), -T ) n_0 (x_0) \nonumber \\
e' (-T , \epsilon_0 ) e(-T , \epsilon_0 ) & = 4  \epsilon_0^2 n'_0 (x_0) \mathcal M (x_0, n_0, \epsilon_0, -T) n_0 (x_0), \label{sq_norm}
\end{align}
where, \[ \mathcal M (x_0, n_0, \epsilon_0, -T) := \left (  \frac{\partial \Phi}{\partial x } (x_0 - \epsilon_0 n_0, -T )  \right )' \frac{\partial \Phi}{\partial x } (x_0 - \epsilon_0 n_0, -T ) \] is a positive definite matrix. The singular value decomposition of the matrix $\mathcal M (x_0, n_0, \epsilon_0, -T)$ would give us,
\begin{align}
& \mathcal M (x_0, n_0, \epsilon_0, -T) =  \nonumber \\
& U' (x_0, n_0, \epsilon_0, -T) ~\Sigma (x_0, n_0, \epsilon_0, -T) ~ U (x_0, n_0, \epsilon_0, -T), \label{svd}
\end{align}
where, $U (x_0, n_0, \epsilon_0, -T)$ is a unitary matrix with $i^{th}$ column $u_i (x_0, n_0, \epsilon_0, T)$, and $\Sigma (x_0, n_0, \epsilon_0, -T) = diag\left (\Lambda_1 (x_0, n_0, \epsilon_0, -T) , \dots, \Lambda_n (x_0, n_0, \epsilon_0, -T)  \right)$, where $\Lambda_1 \ge \Lambda_2 \ge \dots \ge \Lambda_n > 0$. This gives us,
 \begin{align}
  & \mathcal M (x_0, n_0, \epsilon_0, -T) = \nonumber \\
  & \sum_{i=1}^N \Lambda_i (x_0, n_0, \epsilon_0, -T) u_i (x_0, n_0, \epsilon_0, -T)  u_i'  (x_0, n_0, \epsilon_0, -T) . \label{svd_2}
 \end{align}
 Combining \ref{sq_norm} and \ref{svd_2}, 
 \begin{align*} 
 & e' (-T, \epsilon_0) e(-T, \epsilon_0)  =  \\
 & 4  \epsilon^2  \sum_{i=1}^N \Lambda_i (x_0, n_0, \epsilon, -T) n'_0 u_i (x_0, n_0, \epsilon, -T)  u_i'  (x_0, n_0, \epsilon, -T) n_0. 
 \end{align*}
By appropriately selecting  $\epsilon_0$, 
\begin{align*} & e' (-T , \epsilon_0) e(-T , \epsilon_0)  > \\
 & \mathcal {K}(-T, \epsilon_0 ) e' (0 , \epsilon_0) e(0 , \epsilon_0), 
 \end{align*}
where, $\mathcal{K} (-T,\epsilon_0)$ can be made arbitrarily large by making $\epsilon_0$ small.   Also, $e(0, \epsilon_0) = 2 \epsilon n_0, \text {and} ~ e'(0, \epsilon_0) e(0, \epsilon_0) = 4 \epsilon^2 n'_0 n_0 =  4 \epsilon^2 $, as $ n_0' n_0 = 1$. Hence,  
\begin{align*}
 & 4  \epsilon_0^2  \sum_{i=1}^n \Lambda_i (x_0, n_0, \epsilon_0, -T) n'_0 u_i (x_0, n_0, \epsilon_0, -T)  u_i'  (x_0, n_0, \epsilon_0, -T) n_0 \\ 
 & > 4  \epsilon_0^2 \mathcal{K} (-T ,\epsilon_0), \\
 & \sum_{i=1}^n \Lambda_i (x_0, n_0, \epsilon_0, -T) n'_0 u_i (x_0, n_0, \epsilon_0, -T)  u_i'  (x_0, n_0, \epsilon_0, -T) n_0  \\
 & > \mathcal{K} (-T ,\epsilon_0)  , \\
 &  \sum_{i=1}^n \Lambda_i (x_0, n_0, \epsilon_0, -T) \left ( n'_0 u_i (x_0, n_0, \epsilon_0, -T) \right )^2  \\
 & > \mathcal{K} (-T ,\epsilon_0)  . 
\end{align*}
By choosing, $n(x_0) = u_n$, we get, $\Lambda_n (x_0, n_0, \epsilon_0, -T)  > \mathcal{K} (-T ,\epsilon) $. 
For the unit normal vector $\eta_0 (x_0)$, we have, 
\begin{align*} 
& \frac{1}{  \sum_{i=1}^n \frac{1}{\Lambda_i (x_0, n_0, \epsilon_0, -T)} \left ( \eta'_0 u_i (x_0, n_0, \epsilon_0, -T) \right )^{ 2} }  = \\
& \frac{1}{  \sum_{i=1}^n \frac{1}{\Lambda_i (x_0, n_0, \epsilon_0, -T)} \eta'_0 \eta_0 } = \\
&   \frac{1}{  \sum_{i=1}^N \frac{1}{\Lambda_i (x_0, n_0, \epsilon_0, -T)} }  > \frac{\Lambda_n}{n} >  \frac{\mathcal{K} (-T ,\epsilon_0)}{n} .   
\end{align*}  
The $\epsilon_0$ can be made sufficiently small such that, $  \mathcal{K} (- T ,\epsilon_0) > n $. This will give, 
\begin{align}
\frac{1}{  \sum_{i=1}^n \frac{1}{\Lambda_i (x_0, n_0, \epsilon_0, -T)} \left ( \eta'_0 u_i (x_0, n_0, \epsilon_0, -T) \right )^{ 2} } & > 1, \\
 \frac{1}{  \eta_0^T  \left ( \left( \nabla  \Phi(x_0 - \epsilon_0 n_0 , -T)) \right)^T \nabla  \Phi(x_0 - \epsilon_0 n_0 , - T)  \right ) ^{-1} \eta_0 } & > 1, \\
 \rho(x_0 - \epsilon_0 n_0, T) > 1.
\end{align}
Similarly,  we use,
\[ e(- T, \epsilon_0)  =  - 2 \epsilon \frac{\partial \Phi}{\partial x } (x_0 +  \epsilon_0 n_0 (x_0), -T ) n_0 (x_0)  \] and use the successive steps to prove $ \rho(x_0 -  \epsilon_0 n_0, - T) > 1$ . Now let us define a set, \[ \mathcal{N}_{\epsilon} (x^s) := \{ ~ x ~ | ~ \parallel x - x^s \parallel \le \epsilon \} ,\] where $\epsilon < \epsilon_0$.   This gives us, $\rho(x, - T) > 1$ for all $x$, which is at most $\epsilon$ away from the point $x^s$. Hence the proof.
\end{proof}
Below, we provide the proof of Proposition \ref{LE_prop_1}. \\
\begin{proof}
Let us define \[ \zeta(t) := \Phi(x_0, t + \Delta t) - \Phi(x_0, t), \] for initial condition $x_0 \in A(x^s) \setminus x^s$, and $\Delta t $ is a small positive number. It can be observed that $\zeta(0) \neq 0$, as $x_0 \in A(x^s) \setminus x^s$. Also, it can be observed $\Phi( x_0 , \Delta t) \in A(x^s) \setminus x^s$. Considering two initial conditions $x_0$, and $\Phi(x_0, \Delta t)$ the following assertions can be made - there exists a $T_1^*$ such that $\parallel \Phi (x_0, t) - x^s \parallel < \frac{\parallel \zeta(0) \parallel}{2}$ for all $t \ge T_1^*$. Also, there exists a $T_2^*$ such that $\parallel \Phi (x_0, t + \Delta t ) - x^s \parallel < \frac{\parallel \zeta(0) \parallel}{2}$ for all $t \ge T_2^*$. We defining $T^* := max(T_1^*, T_2^*)$.  Using triangle inequality 
\begin{align*} 
& \parallel \Phi (x_0, t + \Delta t ) -  \Phi (x_0, t  ) \parallel \\
& < \parallel \Phi (x_0, t + \Delta t ) -  x^s \parallel + \parallel x^s -  \Phi (x_0, t  ) \parallel < \parallel \zeta_0 \parallel  , ~ \forall ~ t \ge T^*. 
\end{align*}
Hence, we can assert that,  $\parallel \Phi (x_0, t + \Delta t ) -  \Phi (x_0, t  ) \parallel < \parallel \zeta(0) \parallel$ for all $t \ge T^*$. This gives us, $ \lambda (x_0, 0 , t) < 0 $ for all $t \ge T^*$.    
\end{proof}
\bibliographystyle{plain}
\bibliography{ref_new}
\end{document}